\documentclass{amsart}
\title{Preservation of algebraicity in free probability}
\author{Greg W. Anderson}
\address{School of Mathematics, University of Minnesota, Minneapolis, MN 55455}
\email{gwanders@umn.edu}
\date{August 14, 2014}
\newcommand{\Mbold}{{\mathbf{M}}}

\newcommand{\T}{{\mathrm{T}}}

\newcommand{\gbold}{{\mathbf{g}}}

\newcommand{\HHH}{{\mathcal{H}}}

\newcommand{\zero}{{\mathbf{0}}}
\newcommand{\abold}{{\mathbf{a}}}
\newcommand{\bbold}{{\mathbf{b}}}
\newcommand{\cbold}{{\mathbf{c}}}
\newcommand{\dbold}{{\mathbf{d}}}
\newcommand{\val}{{\mathrm{val}}}
\newcommand{\ord}{{\mathrm{ord}}}
\newcommand{\KK}{{\mathbb{K}}}
\newcommand{\NN}{{\mathbb{N}}}
\newcommand{\Xbold}{{\mathbf{X}}}
\newcommand{\GL}{{\mathrm{GL}}}
\newcommand{\Mfrak}{{\mathfrak{M}}}

\newcommand{\AAA}{{\mathcal{A}}}
\newcommand{\ebold}{{\mathbf{e}}}

\newcommand{\ZZ}{{\mathbb{Z}}}
\newcommand{\one}{{\mathbf{1}}}
\newcommand{\trace}{{\mathrm{tr}\,}}
\newcommand{\CC}{{\mathbb{C}}}
\newcommand{\RR}{{\mathbb{R}}}
\newcommand{\Mat}{{\mathrm{Mat}}}
\newtheorem{Proposition}[subsubsection]{Proposition}
\newtheorem{Lemma}[subsubsection]{Lemma}
\newtheorem{Theorem}{Theorem}

\newcommand{\QQ}{{\mathbb{Q}}}

\begin{document} 
\begin{abstract}
We show that any  matrix-polynomial combination of free noncommutative random variables each having an algebraic
law has again an algebraic law. Our result answers a question 
raised by a recent paper of Shlyakhtenko and Skoufranis.
The result belongs to a family of results with origins outside
free probability theory, including a result of Aomoto asserting
algebraicity of the Green function
of random walk of quite general type on a free group.
\end{abstract}
\maketitle
\tableofcontents

\section{Statement of the main result and introduction}

Our main result is as follows:

\begin{Theorem}\label{Theorem:MainResult}Let $(\AAA,\phi)$ be a noncommutative 
probability space. Let
$x_1,\dots,x_q\in \AAA$
be freely independent noncommutative random variables. Let
$$X\in \Mat_p(\CC\langle x_1,\dots,x_q\rangle)\subset \Mat_p(\AAA)$$
be a matrix. If the laws of $x_1,\dots,x_q$ are algebraic, then so is the law of $X$.
\end{Theorem}

This paper is devoted to a proof of Theorem \ref{Theorem:MainResult}. 
We say for short that $X$ is a {\em free matrix-polynomial combination} of $x_1,\dots,x_q$. 
See \S\ref{section:Breakdown} below for notation, terminology, background and references clarifying the theorem statement.

The beginning student of free probability immediately notices that all the 
important special distributions in the theory are algebraic, e.g., 
the arcsine law, Wigner's semicircle law, the Marcenko-Pastur law, and so on.
Our main result and its proof describes a mechanism at least in part accounting for this ubiquity.

The phenomenon of preservation of algebraicity 
established in general here has previously been verified  
in many special cases, often merely as a byproduct.
Here are several important examples.
(i) It is implicit in the theory of the $R$-transform introduced in \cite{Voiculescu} that free (additive) convolution preserves algebraicity. 
We note that this phenomenon has been exploited in a practical way in \cite{EdelmanRao}.  Similar remarks apply to free multiplicative convolution.
(ii) It is implicit in
the theory of commutators of free random variables developed in \cite{NicaSpeicher}
that formation of free commutators preserves algebraicity. 
(iii) Algebraicity of the law of a free matrix-polynomial combination of semicircular variables is asserted in \cite[Thm. 5.4]{ShlyakhSkou}.
This result forms part of a result of wide scope, namely \cite[Thm. 1.1]{ShlyakhSkou}, which gives important constraints
on the structure of the law of a free matrix-polynomial combination of noncommutative random variables each of which has a nonatomic law.

Our main result answers a question raised by \cite{ShlyakhSkou}.

One further example of algebraicity, originating outside free probability, deserves special mention as the archetype for Theorem \ref{Theorem:MainResult}.

\begin{Theorem}[See, e.g., {\cite[Cor. 6.7.2, p. 210]{Stanley}}]\label{Theorem:Aomoto}
Let $G$ be a free group.
Let $\CC[G]$ be the group algebra of $G$ with complex coefficients. 
Let $\tau:\CC[G]\rightarrow\CC$ be the unique $\CC$-linear map such that $\tau(g)=\delta_{1g}$ for $g\in G$.
Then for any $P\in \CC[G]$ and $g\in G$ the formal power series $\sum_{n=0}^\infty \tau(gP^n)t^n\in \CC[[t]]$ is algebraic over the field of rational functions $\CC(t)$.
\end{Theorem} 

It is remarked in \cite{Kontsevich} that Theorem \ref{Theorem:Aomoto} has been frequently rediscovered.   We do not know the identity of the first discoverer.

Now a unitary noncommutative random variable factors as a product of two free Bernoulli variables.
In other plainer words, a group generated by elements $y_1,\dots,y_{2q}$ subject only to the relations
$y_i^2=1$ for $i=1,\dots,2q$ contains a free subgroup on $q$ generators,
for example that generated by $y_1y_2,y_3y_4,\dots,y_{2q-1}y_{2q}$.
Thus Theorem \ref{Theorem:Aomoto} in the case $g=1$ is 
a consequence of Theorem \ref{Theorem:MainResult}. 
Given this close relationship, 
one reasonably looks to the proof of Theorem \ref{Theorem:Aomoto} for 
clues concerning the proof of Theorem \ref{Theorem:MainResult}.

The proof of Theorem \ref{Theorem:Aomoto} given in \cite{Stanley} is based on \cite{Haiman} which in turn builds upon the theory of algebraic 
noncommutative formal power series (see e.g., \cite{SalSoi} or
\linebreak \cite[Chap. 6]{Stanley}), and for the latter, crucial foundations are laid in the seminal paper \cite{ChomSchu} on context free languages.
Tools from the same kit are also used to prove \linebreak \cite[Thm. 5.4]{ShlyakhSkou}. 

Somewhat counterintuitively, we prove Theorem \ref{Theorem:MainResult} by an approach mostly avoiding formal language theory.
We rely instead on methods from free probability, random walk on groups, algebraic geometry and commutative algebra. 
Formal language theory is still involved, but in a different and simpler way. A proof of Theorem \ref{Theorem:MainResult}
parallel to that of \cite[Thm. 5.4]{ShlyakhSkou} might yet be possible and would be very interesting.
For now the sticking point seems to be that no suitable generalization of \linebreak \cite[Lemma 5.12]{ShlyakhSkou}
is obviously on offer. It is conceivable that our methods could have in the reverse direction some
impact on formal language theory.

Questions quite similar to that answered by Theorem \ref{Theorem:Aomoto} 
have been treated in the literature of random walk on groups.
Here are two particularly important examples.
(i) In \cite{Aomoto} algebraicity of the Green function of random walk
of a fairly general type on a finitely generated free group  was proved by 
explicit calculation.
This is the earliest paper of which we are aware which 
up to some mild and ultimately removable hypotheses proves Theorem \ref{Theorem:Aomoto}.
(ii) In \cite{Woess} algebraicity of the Green function of any finitely supported random walk
on a group with a finitely generated free subgroup of finite index  was proved
by a method based on formal language theory, a method much the same as later used in \cite{Stanley} 
to prove Theorem \ref{Theorem:Aomoto}.

The overall approach to algebraicity used here is 
adapted not from the literature of formal language theory,
but rather from that of random walk on groups and graphs, especially free groups and infinite trees.
 See \cite{WoessBook} for background.
From the random walk literature we cite as particular examples \cite{Aomoto}, \cite{Lalley1} and \cite{NagWoess}.
This short list could be greatly extended because the idea we are adapting is a fundamental trope.
There are two main components to this trope.
Firstly, one  gets relatively simple recursions for Green functions by exploiting finiteness of cone types, or some related principle of self-similarity. Such recursions are often (in effect)  consequences of the familiar matrix inversion formula
$$
\left[\begin{array}{cc}
\abold&\bbold\\
\cbold&\dbold\end{array}\right]^{-1}=
\left[\begin{array}{cc}
\zero&\zero\\
\zero&\dbold^{-1}
\end{array}\right]+\left[\begin{array}{r}
\one\\
-\dbold^{-1}\cbold
\end{array}\right]
(\abold-\bbold \dbold^{-1}\cbold)^{-1}\left[\begin{array}{rr}
\one&-\bbold \dbold^{-1}
\end{array}\right].
$$
Here we put the Boltzmann-Fock space model of free random variables \cite{Voiculescu}
into suitably ``arboreal'' form in order to gain access to similar recursions.
Secondly, there are criteria available for recognizing when recursions 
have algebraic solutions. Of course formal language theory provides such criteria, but
 there are also less complicated criteria. For example, \cite[Prop. 5.1]{Lalley1} is an especially clear and general criterion,
and it is applied in the cited paper in an instructive manner.
We use a similar but simpler criterion here. See Proposition \ref{Proposition:AlgebraicityCriterion} below.

We remark that  when writing \cite{AndersonZeitouni}, wherein was presented an algebraicity criterion
similar to if rather more complicated than \cite[Prop. 5.1]{Lalley1},
namely \cite[Thm. 6.1]{AndersonZeitouni},
the authors were unfortunately unaware of \cite{Lalley1}. We wish to acknowledge the priority.

We also remark that relations between algebraicity and positivity are highly developed in the random walk literature,
e.g., in \cite{Lalley1} and \cite{NagWoess}, 
leading to local limit theorems. We do not touch those ideas here but we think they could be fruitfully applied in the free probability context.

The paper \cite{BelMaiSpei} has been an important influence
because, building upon operator-valued $R$-transform theory \cite{VoiculescuBis}, \cite{Dykema}
it reveals a rich algebraic and analytic structure to exploit for studying algebraicity.
In particular, study of the fixed point equations
stated in \cite[Thm. 2.2]{BelMaiSpei} should in principle lead to a proof of Theorem \ref{Theorem:MainResult}.
But because the many-variable setting for these equations is too difficult for us to handle, we work instead with a lightly modified version of the original setup of \cite{Voiculescu} and we use just the single classical parameter $z$. 
It remains an open and interesting problem to prove algebraicity in the many-variable setting
of \cite{BelMaiSpei}.

The linearization trick which we learned from the papers \cite{HST} and \cite{HT},
which we refined and used in \cite{Anderson}, and which in the refined form was also used in \cite{BelMaiSpei} plays here an important role as well. But as we have recently learned 
from \cite{HMV06}
and wish to acknowledge here,
the  trick in its refined form already exists in the literature.
A theorem of Sch\"{u}tzenberger \cite{Schutz} (see also \cite[Thm. 7.1, p. 18]{BerReu})  belonging to  formal language theory
contains the core idea of the trick sans the self-adjointness-preserving aspect. The whole trick is contained in \cite[Lemma 4.1]{HMV06}  and is called there {\em symmetric realization}.  (There is a minor issue that the lemma is stated in the real case but the hermitian generalization is routine.)
The cited lemma appears in a context perhaps not at first glance closely connected with free probability,
but clearly and closely allied with linear systems theory 
\cite{Brockett}, \cite{Kalman}. The author plans with several co-authors 
to report on these interconnections in a forthcoming paper.

As it happens, because Theorem \ref{Theorem:MainResult} does not mention self-adjointness, all we  need 
of the linearization/realization technique is Sch\"{u}tzenberger's theorem and of the latter we in fact need only a  fragment, which we prove here ``from scratch'' in a paragraph. See Lemma \ref{Lemma:Linearization} below.
See also \S\ref{subsubsection:Automaton} for an explanation of the connection with Sch\"{u}tzenberger's theorem.
Thus  formal language theory does  enter into our proof of Theorem \ref{Theorem:MainResult}
 after all, but in a rather simple way.

In the paper \cite{KassReut}, generalizing a result of the earlier paper \cite{Kontsevich},
the following variant of Theorem \ref{Theorem:Aomoto} is proved.
In this variant, instead of a group-ring element
\linebreak $P\in \CC[G]$ with complex coefficients,
one considers a square matrix 
$P\in \Mat_n(\ZZ[G])$ of group-ring elements with integer coefficients,
one considers not the usual generating function
but rather the {\em zeta-function}
$$\exp\left(\sum_{k=1}^\infty \tau(\trace(P^k))\frac{t^k}{k}\right)\in 1+t\CC[[t]],$$
and finally, one concludes strikingly that the latter both has integer coefficients
and is algebraic. Thus motivated, we raise the following question.
Is ``integral-algebraicity of zeta-functions''
 preserved under ``integral-matrix-polynomial combination'' of 
free random variables? Perhaps this question could be answered by combining methods used here with those of the cited papers.

As mentioned above, our formal algebraic setup for proving Theorem \ref{Theorem:MainResult}
is based on the original setup of \cite{Voiculescu}. From that starting point,
we make the rest of our definitions so as to keep our approach
to proving Theorem \ref{Theorem:MainResult} as simple as possible. 
We make no positivity assumptions---moment sequences of variables can be arbitrarily prescribed sequences of complex numbers. 
We work over the field $\CC((1/z))$ of formal Laurent series, using simple ideas about Banach algebras over complete ultrametrically normed fields in lieu of operator theory over the complex numbers. Our method reveals nothing about the branch points of the algebraic functions it produces. 
It is an open problem to recover information about positivity and branch points.
Perhaps this is only a matter of unifying features of the several theories mentioned above, but in our opinion some further ingredients from algebraic geometry 
will be needed. The soliton theory literature, e.g., \cite{Mumford}, might provide guidance.

Here is an outline of the paper. In \S\ref{section:Breakdown}, after filling in background in leisurely fashion,
and in particular writing down a simple algebraicity criterion,
namely Proposition \ref{Proposition:AlgebraicityCriterion} below,
we reformulate Theorem \ref{Theorem:MainResult} as the conjunction of two propositions both of which concern the generalized Schwinger-Dyson equation.
In \S\ref{section:Gentle} we introduce the formal algebraic variant of operator theory used in this paper. In \S\ref{section:Linearization} we introduce a suitable model of free random variables based on that of \cite{Voiculescu}  and  we deploy the linearization/realization technique, a.k.a. Sch\"{u}tzenberger's theorem.
In \S\ref{section:GSDE} we exhibit the solutions of the generalized Schwinger-Dyson equation needed to prove Theorem \ref{Theorem:MainResult}.
In the remainder of the paper we switch to the viewpoint of
algebraic geometry and commutative algebra. In \S\ref{section:NewtonPolygon} we review 
topics connected with singularities of plane algebraic curves, especially Newton polygons. In \S\ref{section:Gizmo}
we apply the Weierstrass Preparation Theorem in a perhaps unexpected way.
Finally, in \S\ref{section:Endgame} we complete the proof of Theorem \ref{Theorem:MainResult}
by checking hypotheses in Proposition \ref{Proposition:AlgebraicityCriterion}.

Lastly, we remark that the paper is long only because we have included many explanations and 
reviews of background to smooth the way for the interested reader who might not be familiar with all the (seemingly) disparate materials collected here.

\section{Background for the main result and a reduction of the proof}
\label{section:Breakdown}
After recalling principal definitions, fixing notation, 
and filling in background for Theorem \ref{Theorem:MainResult},  we reduce Theorem \ref{Theorem:MainResult}  to two propositions each treating some aspect 
of the generalized Schwinger-Dyson equation.

\subsection{Noncommutative probability spaces and free independence}
We \linebreak present a brief review to fix notation.
See, e.g., \cite{AGZ}, \cite{NicaSpeicherBis}, \cite{SpeicherICM}, or \cite{VDN} for background.

\subsubsection{Algebras}
All algebras in this paper are unital, associative,
and have a scalar field containing $\CC$. 
The unit of an algebra $\AAA$ is denoted by $1_\AAA$;
other notation, e.g., simply $1$,  may be used when context permits.
Given elements $x_1,\dots,x_q\in \AAA$ of an algebra,
let $\CC\langle x_1,\dots,x_q\rangle
\subset \AAA$ denote the subring of $\AAA$ 
generated by forming all finite $\CC$-linear combinations
of monomials in the given elements $x_1,\dots,x_q$, including
the ``empty monomial'' $1_\AAA$. (But if $\AAA$ is commutative,
instead of $\CC\langle x_1,\dots,x_q\rangle$, we write $\CC[x_1,\dots,x_q]$ 
as is usual in commutative algebra.)
 
\subsubsection{Noncommutative probability spaces}
A {\em state} $\phi$ on an algebra $\AAA$ is simply a $\CC$-linear functional $\phi:\AAA\rightarrow\CC$
such that $\phi(1_\AAA)=1$. In our formal algebraic setup no positivity constraints are imposed. 
A {\em noncommutative probability space} is a pair $(\AAA,\phi)$ consisting of an algebra $\AAA$
and a state $\phi$ on that algebra.  Given such a pair $(\AAA,\phi)$,  elements of $\AAA$
are called {\em noncommutative random variables}. 

\subsubsection{Matrices with algebra entries}
Given an algebra $\AAA$ and a positive integer $n$, let $\Mat_n(\AAA)$ denote the  algebra of
$n$-by-$n$ matrices with entries in $\AAA$.
 More generally let $\Mat_{k\times \ell}(\AAA)$ denote the space of $k$-by-$\ell$ matrices with entries in 
$\AAA$. For \linebreak $A\in \Mat_{k\times \ell}(\CC)$ and $a\in \AAA$ we define
$A\otimes a\in \Mat_{k\times \ell}(\AAA)$ by $(A\otimes a)(i,j)=A(i,j)a$.
Let $1=I_n=I_n\otimes 1_\AAA\in \Mat_n(\AAA)$ denote the
identity matrix, as context may permit.
Let  $\ebold_{ij}\in \Mat_{k\times \ell}(\CC)$ 
denote the elementary matrix with $1$ in position $(i,j)$ and $0$ in every other position.
 Let $\GL_n(\AAA)$ denote the group of invertible
elements of $\Mat_n(\AAA)$.  
Given a noncommutative probability
space $(\AAA,\phi)$, we regard each matrix $A\in \Mat_n(\AAA)$
as a noncommutative random variable with respect to the state  $\phi_n:\Mat_n(\AAA)\rightarrow\CC$
given by the formula $\phi_n(A)=\frac{1}{n}\sum_{i=1}^n \phi(A(i,i))$.

\subsubsection{Free independence}
Let $(\AAA,\phi)$ be a noncommutative probability space and let $\AAA_1,\dots,\AAA_q\subset \AAA$ be subalgebras such that $1_\AAA\in \cap_{i=1}^q \AAA_i$.
One says that $\AAA_1,\dots,\AAA_q$ are {\em freely independent} if for every positive integer $k$,
sequence $i_1,\dots,i_k\in \{1,\dots,q\}$ such that 
$i_1\neq i_2$, $i_2\neq i_3$, \dots, $i_{k-1}\neq i_k$
and sequence $x_1\in \AAA_{i_1}$, \dots, $x_k\in \AAA_{i_k}$ such that
$\phi(x_1)=\cdots=\phi(x_k)=0$,
one has $\phi(x_1\cdots x_k)=0$.
As a special case of the preceding general definition, one says that noncommutative random variables $x_1,\dots,x_q\in \AAA$ are {\em freely independent} if
the subalgebras $\CC\langle x_1\rangle,\dots,\CC\langle x_q\rangle\subset \AAA$ are freely independent.

\subsubsection{Univariate laws}
Let $\Xbold$ be a variable. A {\em univariate law} (or, context permitting, simply a {\em law}) is by definition a 
state $\mu:\CC\langle\Xbold\rangle\rightarrow\CC$ on the one-variable polynomial algebra $\CC\langle\Xbold\rangle$. 
The value $\mu(\Xbold^n)\in \CC$ is called the {\em $n^{th}$ moment} of $\mu$.
Note that in our formal algebraic setup the moments of a law are allowed to be arbitrarily prescribed complex numbers.
Given a noncommutative probability space $(\AAA,\phi)$ and a noncommutative random variable $x\in \AAA$,
the {\em law} of $x$ is by definition the linear functional $\mu_x:\CC\langle\Xbold\rangle\rightarrow\CC$
determined by the formula $\mu_x(\Xbold^n)=\phi(x^n)$ for integers $n\geq 0$.

\subsubsection{Noncommutative joint laws}\label{subsubsection:RemarkNC}
Let $\Xbold_1,\dots,\Xbold_q$ be independent noncommuting algebraic variables
and let $\CC\langle \Xbold_1,\dots,\Xbold_q\rangle$ be the noncommutative polynomial ring generated by these variables. 
A {\em $q$-variable noncommutative law} or, context permitting, simply a {\em law}, is a state on the algebra $\CC\langle \Xbold_1,\dots,\Xbold_q\rangle$.
Let $(\AAA,\phi)$ be a noncommutative probability space and let $x_1,\dots,x_q\in \AAA$
be noncommutative random variables. 
The {\em joint law}
$\mu_{x_1,\dots,x_q}:\CC\langle \Xbold_1,\dots,\Xbold_q\rangle\rightarrow\CC$ of the $q$-tuple $(x_1,\dots,x_q)$ is by definition the linear functional defined by the  rule
$\mu_{x_1,\dots,x_q}(f(\Xbold_1,\dots,\Xbold_q))=\phi(f(x_1,\dots,x_q))$ for $f(\Xbold_1,\dots,\Xbold_q)\in \CC\langle \Xbold_1,\dots,\Xbold_q\rangle$.
The laws $\mu_{x_1},\dots,\mu_{x_q}$ of the individual variables (by analogy with classical probabilistic usage) are 
 called the {\em marginal} laws for the joint law $\mu_{x_1,\dots,x_q}$.
A point worth emphasizing is that if
$x_1,\dots,x_q$ are freely independent, then the joint law $\mu_{x_1,\dots,x_q}$
 is uniquely determined by the marginal laws $\mu_{x_1},\dots,\mu_{x_q}$.

\subsection{The Laurent series field $\CC((1/z))$ and related notions} 
We recall several definitions together providing a framework in which to discuss algebraicity.
See the text \cite{Artin} by Artin for background on valued fields and algebraic functions.

\subsubsection{Definition of $\CC((1/z))$ and related objects}
Let $\CC((1/z))$ denote the set of {\em formal Laurent series} in $z$ of the form
\begin{equation}\label{equation:fLaurent}
f=\sum_{i\in \ZZ}c_i z^i\;\;
(\mbox{$c_i\in \CC$ and $c_i=0$ for $i\gg 0$}).
\end{equation}
(The coefficients $c_i$ are not subject to any majorization.)
Equipped with addition and multiplication in evident fashion, the set $\CC((1/z))$ becomes a field.
Note that we have inclusions 
$$\CC[z]\subset \CC(z)\subset \CC((1/z))\;\mbox{and}\;\CC[[1/z]]\subset \CC((1/z))$$
where $\CC[z]$ is the ring of polynomials in $z$, $\CC[[1/z]]$ is the ring of formal power series in $1/z$,
and $\CC(z)$ is the field of rational functions of $z$, all with coefficients in $\CC$.
Note also that we have an additive direct sum decomposition
$$\CC((1/z))=\CC[z]\oplus (1/z)\CC[[1/z]].$$
In our algebraic setup the formal variable $z$ corresponds to the classical parameter $z$ in the upper half-plane.

\subsubsection{Algebraic elements of $\CC((1/z))$ and their irreducible equations}\label{subsubsection:Algebraic}
Let $\CC[x,y]$
be the polynomial ring over $\CC$ in two independent (commuting) variables $x$ and $y$.
We say that $f\in \CC((1/z))$ is {\em algebraic} if one and hence all three of the following equivalent conditions
hold:
\begin{itemize}
\item There exists some $0\neq P(x,y)\in \CC[x,y]$
such that $P(z,f)=0$.
\item There exists some $0\neq Q(x,y)\in \CC[x,y]$
such that $Q(1/z,f)=0$.
\item The field $\CC(z,f)$ generated over $\CC(z)$ by $f$ is a vector space of finite dimension over $\CC(z)$. 
\end{itemize}
As is well-known,
the algebraic elements form a subfield of $\CC((1/z))$ containing $\CC(z)$.
For algebraic $f\in \CC((1/z))$ there exists irreducible $F(x,y)\in \CC[x,y]$
unique up to a constant multiple such that $F(1/z,f)=0$. (The insertion of $1/z$ in the preceding definition
rather than $z$ is a technical convenience.)
With but slight abuse of language
we call any such irreducible polynomial the {\em irreducible equation} of $f$.

\subsubsection{Valuations}
For $f\in \CC((1/z))$ expanded as on line \eqref{equation:fLaurent} 
we define
$$\val\,f =\sup\, \{i\in \ZZ\mid c_i\neq 0\}=(\mbox{the {\em valuation} of $f$})\in \ZZ\cup\{-\infty\}.$$
Note that 
\begin{eqnarray}\label{equation:VP1}
\val\,f=-\infty&\Leftrightarrow& f=0,\\\label{equation:VP2}
\val(f_1f_2)&=&\val\,f_1+\val\,f_2,\\\label{equation:VP3}
\val(f_1+f_2)&\leq&\max(\val\, f_1,\val\, f_2)\;\;\mbox{with equality if $\val\,f_1\neq \val\,f_2$.}
\end{eqnarray}
Thus $\val$ is (the logarithm of) a nonarchimedean valuation in the sense of \cite{Artin}.
Thus it becomes possible to use (ultra)metric space ideas to reason about $\CC((1/z))$ and related objects, as in \cite{Artin}, and we will do so throughout this paper. We may speak for example of completeness.
It is easy to see that $\CC((1/z))$ is complete with respect to the valuation $\val$.

\subsubsection{Banach algebra structure for $\Mat_n(\CC((1/z)))$}
We now extend the metric space ideas a bit farther. We equip the matrix algebra $\Mat_n(\CC((1/z)))$
with a valuation by the rule $\val\, A=\max_{i,j=1}^n \val\,A(i,j)$. 
Then \eqref{equation:VP1} and \eqref{equation:VP3} hold for matrices,
\eqref{equation:VP2} holds for multiplication of a matrix by a scalar,
and \eqref{equation:VP2} holds for multiplication
of two matrices provided that ``$=$'' is relaxed to ``$\leq $.''
Thus $\Mat_n(\CC((1/z)))$ becomes a Banach algebra over $\CC((1/z))$.
Later, in \S\ref{section:Gentle}, a certain infinite-dimensional Banach  algebra over $\CC((1/z))$ 
will be introduced.

\subsubsection{Composition of Laurent series}\label{subsubsection:Composition}
The composition $f\circ g\in \CC((1/z))$
of\linebreak $f,g\in \CC((1/z))$ is defined provided that $\val\, g>0$.
The set $z+\CC[[1/z]]$ forms a group under composition.
This group acts on the right side of $\CC((1/z))$ by $\CC$-linear field automorphisms.

\begin{Lemma}\label{Lemma:InversionAndAlgebraicity}
If $f,g\in z+\CC[[1/z]]$ satisfy $f\circ g=z$ and $f$ is algebraic,
then $g$ is also algebraic.
\end{Lemma}
\proof
For some $0\neq P(x,y)\in \CC[x,y]$ we have  
$0=P(z,f)\circ g=P(g,z)$.
\qed

\subsection{Algebraicity of univariate laws}
We recall how to attach to each univariate law a (formal) Stieltjes transform and a  (modified formal) $R$-transform {\em \`{a} la} Voiculescu. 
Then we recall how in terms of these transforms one can characterize algebraicity of a law.
\subsubsection{Formal Stieltjes transforms}
For a law $\mu:\CC\langle\Xbold\rangle\rightarrow\CC$,
the formal sum
$$S_\mu(z)=\sum_{n=0}^\infty \mu(\Xbold^n)/z^{n+1}\in \CC((1/z))$$ 
is by definition the {\em formal Stieltjes transform} of $\mu$.
Hereafter we drop the adjective ``formal'' since no other kind of Stieltjes transform will be considered in this paper.

\subsubsection{Algebraicity of univariate laws}
A law $\mu$  will be called {\em algebraic} if its  Stieltjes transform $S_\mu(z)\in \CC((1/z))$ is algebraic.

\subsubsection{Free cumulants and $R$-transforms}
Given a law $\mu:\CC\langle\Xbold\rangle\rightarrow\CC$ one defines in free probability theory for each positive integer $n$ the {\em $n^{th}$ free cumulant} $\kappa_n(\mu)\in \CC$. This can be done  various ways, e.g., with generating functions 
or combinatorially using noncrossing partitions. See, e.g.,  \cite{AGZ}, \cite{NicaSpeicherBis}, or \cite{VDN} for background; the founding document for this theory is \cite{Voiculescu}.
The generating function
$$R_\mu(t)=\sum_{n=1}^\infty \kappa_n(\mu)t^{n-1}\in \CC[[t]]$$
for the free cumulants is the formal version of the {\em $R$-transform} of Voiculescu. 
Hereafter we drop the adjective ``formal'' since no other kind of $R$-transform will be considered in this paper.

\subsubsection{Modified $R$-transforms}
To define and make calculations with free cumulants, we will use the generating function method.
Consider the {\em modified} $R$-transform
 $$\tilde{R}_\mu(z)=z+R_\mu(1/z)=z+\sum_{n=1}^\infty \kappa_n(\mu)z^{1-n}\in z+\CC[[1/z]],$$
 which we will find slightly more convenient. Obviously each of $R_\mu(t)$ and $\tilde{R}_\mu(z)$
uniquely determines the other.
It is known (see \cite{AGZ}, \cite{NicaSpeicherBis}, \cite{Voiculescu}, or \cite{VDN}) that $\tilde{R}_\mu(z)$ is the unique solution of the equation
\begin{equation}\label{equation:RGenFunDef}
\left(\frac{1}{S_\mu(z)}\right)\circ\widetilde{R}_\mu(z) =z.\;\;
\left(\mbox{Equivalently:}\;\;\widetilde{R}_\mu(z)\circ \left(\frac{1}{S_\mu(z)}\right)=z.\right)
\end{equation}
Since $z+\CC[[1/z]]$ is a group under composition, the modified  $R$-transform $\widetilde{R}_\mu(z)$
is well-defined for every law $\mu$, hence the sequence $\{\kappa_n(\mu)\}_{n=1}^\infty$ of free cumulants is defined,
and it uniquely determines $\mu$.
Furthermore the free cumulants of a law can be arbitrarily prescribed.

The next lemma expresses algebraicity in terms of free cumulants.
\begin{Lemma}\label{Lemma:Centering}
Let $\mu:\CC\langle\Xbold\rangle\rightarrow\CC$ be a law.
Then the following statements are equivalent:
\begin{enumerate}
\item[(I)] $\mu$ is algebraic.
\item[(II)] $S_\mu(z)$ is algebraic.
\item[(III)] $\widetilde{R}_\mu(z)$ is algebraic.
\end{enumerate}
\end{Lemma}
\proof The equivalence (I)$\Leftrightarrow$(II) holds by definition.
The equivalence (II)$\Leftrightarrow$(III) holds by 
Lemma \ref{Lemma:InversionAndAlgebraicity}, statement \eqref{equation:RGenFunDef} above,
and the fact that $0\neq f\in \CC((1/z))$ is algebraic if and only if the reciprocal $1/f$ is algebraic.
\qed

\subsection{An algebraicity criterion}\label{section:Algebraicity}
We now present the algebraicity criterion which we will use to take the final step of the proof
of Theorem \ref{Theorem:MainResult}. To do so we abruptly switch to the optic of commutative algebra.
\subsubsection{Setup for the  criterion}
Let $K/K_0$ be any extension of fields.  Let  \linebreak
$x=(x_1,\dots,x_n)$ be an $n$-tuple of independent (commuting) variables
and let $K_0[x]$ be the polynomial ring generated over $K_0$ by these variables.
Let 
$f=(f_1,\dots,f_n)=f(x)\in K_0[x]^n$ be an $n$-tuple of polynomials. 
Let $J(x)=\det_{i,j=1}^n \frac{\partial f_i}{\partial x_j}\in K_0[x]$
be the determinant of the Jacobian matrix of $f$. Let
$\alpha=(\alpha_1,\dots,\alpha_n)\in K^n$
be an $n$-tuple such that $f(\alpha)=0$
but $J(\alpha)\neq 0$.

\begin{Proposition}\label{Proposition:AlgebraicityCriterion}
 Notation and assumptions are as above.  Every entry of the vector
$\alpha$ is algebraic over $K_0$.
\end{Proposition}
This statement is the same as \cite[Prop. 5.3, Chap. VIII, p. 371]{Lang},
and also the same as \cite[Lemma 6.6.9, Chap. 6, p. 206]{Stanley}.

\subsection{Large-scale organization of the proof of Theorem \ref{Theorem:MainResult}}
We recall the generalized Schwinger-Dyson equation and then we state two technical results about it which together imply Theorem \ref{Theorem:MainResult}.
\subsubsection{The generalized Schwinger-Dyson equation}
In the two technical propositions to be formulated below we consider an instance
\begin{equation}\label{equation:GSDE2}
I_{n}+a^{(0)} g+\sum_{\theta=1}^q \sum_{j=2}^\infty \kappa_j^{(\theta)}(a^{(\theta)} g)^j=0
\end{equation}
 of the {\em generalized Schwinger-Dyson equation} for which the data are
\begin{equation}\label{equation:Data} 
\left\{\begin{array}{l}
\mbox{positive integers $q$ and $n$,}\\
\mbox{matrices $g,a^{(0)}\in \Mat_n(\CC((1/z)))$,}\\
\mbox{matrices $a^{(1)},\dots,a^{(q)}\in \Mat_n(\CC)$, and}\\
\mbox{a family $\left\{\left\{\kappa_j^{(\theta)}\right\}_{j=2}^\infty\right\}_{\theta=1}^q$ of complex numbers.}
\end{array}\right.
\end{equation}
We assume that
\begin{equation}\label{equation:GSDE1}
\lim_{j\rightarrow\infty}\val\,(a^{(\theta)} g)^j=-\infty\;\;\mbox{for $\theta=1,\dots,q$}
\end{equation}
in order that the left side of \eqref{equation:GSDE2} have a well-defined value in $\Mat_n(\CC((1/z)))$.
Furthermore, we impose the following nondegeneracy condition: 
\begin{eqnarray}
\label{equation:GSDE3}
&&\mbox{The linear map}\\
\nonumber&&\left(h\mapsto a^{(0)} h+\sum_{\theta=1}^q\sum_{j=2}^\infty
\sum_{\nu=0}^{j-1}
\kappa_j^{(\theta)}
(a^{(\theta)}g)^\nu (a^{(\theta)}h)(a^{(\theta)}g)^{j-1-\nu}
\right)\\
\nonumber&&:\Mat_n(\CC((1/z)))\rightarrow \Mat_n(\CC((1/z)))
\;\mbox{is invertible.}
\end{eqnarray}
Note that the map above is well-defined by assumption \eqref{equation:GSDE1}. 

We will prove the following two results.
\begin{Proposition}\label{Proposition:MainResultBlockManipulation}
Let $(\AAA,\phi)$ be a noncommutative probability space. Let \linebreak
$x_1,\dots,x_q\in \AAA$
be freely independent noncommutative random variables.
Let
$$X\in \Mat_p(\CC\langle x_1,\dots,x_q\rangle)\subset \Mat_p(\AAA)$$
be a matrix. 
(Up to this point we have merely repeated the setup for Theorem \ref{Theorem:MainResult}.) 
 For indices $\theta=1,\dots,q$ and  $j=2,3,4,\dots$, let
$\kappa_j^{(\theta)}$ denote the $j^{th}$ free cumulant of the law of the noncommutative random variable $x_\theta$. Then the family $\left\{\left\{\kappa_j^{(\theta)}\right\}_{\theta=1}^q\right\}_{j=2}^\infty$ of complex numbers
for some integer $n> p$  can be completed to a family 
$$\left(q,n,g,a^{(0)},\left\{a^{(\theta)}\right\}_{\theta=1}^q,\left\{\left\{\kappa_j^{(\theta)}\right\}_{j=2}^\infty\right\}_{\theta=1}^q\right)$$
of the form \eqref{equation:Data}
satisfying \eqref{equation:GSDE2}, \eqref{equation:GSDE1},  and \eqref{equation:GSDE3}
along with the further conditions
\begin{eqnarray}
\label{equation:GSDE5}
&&a^{(0)}\in \{A+Bz\mid A,B\in \Mat_n(\CC)\}\;\;\mbox{and}\\
\label{equation:GSDE0}
&&S_{\mu_X}=-\frac{1}{p}\sum_{i=1}^p g(i,i).
\end{eqnarray}
\end{Proposition}

\begin{Proposition}\label{Proposition:MainResultAlgebraicityProof}
Let data of the form \eqref{equation:Data} satisfy \eqref{equation:GSDE2}, \eqref{equation:GSDE1}, 
and \eqref{equation:GSDE3}.  Assume furthermore that 
\begin{eqnarray}
\label{equation:GSDE5.5}
&&a^{(0)}\in \Mat_n(\CC(z))\;\;\mbox{and}\\
\label{equation:GSDE4}
&&\mbox{$\sum_{j=2}^\infty \kappa_{j+1}^{(\theta)}z^{-j}\in \CC((1/z))$ is algebraic for $\theta=1,\dots,q$.}
\end{eqnarray}
Then every entry of the matrix $g$ is algebraic.
\end{Proposition}
See \S\ref{section:Linearization}  and \S\ref{section:GSDE} below for the proof of Proposition \ref{Proposition:MainResultBlockManipulation}.
See \S\ref{section:NewtonPolygon}, \S\ref{section:Gizmo},
and \S\ref{section:Endgame} below for the proof of Proposition \ref{Proposition:MainResultAlgebraicityProof}.

\subsubsection{Reduction of the proof of Theorem \ref{Theorem:MainResult}}
In view of Lemma \ref{Lemma:Centering}, it is clear that Propositions \ref{Proposition:MainResultBlockManipulation}
and \ref{Proposition:MainResultAlgebraicityProof}
 together imply Theorem \ref{Theorem:MainResult}.

\subsubsection{Remark}
In the simple case 
$$X=x_1+\cdots+x_q\in \CC\langle x_1,\dots,x_q\rangle\subset 
\Mat_1(\AAA)=\AAA,$$ 
the instance of generalized Schwinger-Dyson equation
emerging from the proof of Proposition \ref{Proposition:MainResultBlockManipulation} reduces to the standard fact \cite{Voiculescu} 
that the $R$-transform is additive for the addition of free random variables.
\subsubsection{Remark}
The (un)generalized Schwinger-Dyson equation is familiar in the case that $\kappa_j^{(\theta)}=\kappa_j(\mu_{x_\theta})=0$ for $j>2$.
In the latter special case the equation \eqref{equation:GSDE2}  arises naturally in the study of free matrix-polynomial combinations of semicircular variables. See, e.g., \cite{Anderson}, \cite{AGZ}, \cite{BelMaiSpei}, \cite{HST}, \cite{HT}, \cite{HRS07}, \cite {NicaSpeicherBis} and \cite{VDN}.

\subsubsection{Remark}
Proposition \ref{Proposition:MainResultBlockManipulation} is unsurprising.
It is proved here by straightforwardly
combining three standard methodologies, namely:
\begin{enumerate}
\item[(i)] the Boltzmann-Fock space model of free random variables,
 \item[(ii)] the linearization/realization method, and
 \item[(iii)] recursions of a type occurring
in the study of  random walk on  infinite trees.
\end{enumerate}
Methodology (i) clearly originates in \cite{Voiculescu}.
Methodology (ii) we learned from the papers \cite{HST} and \cite{HT},
and subsequently we refined it in \cite{Anderson}, but as we have recently learned from \cite{HMV06}, the essential point
apart from the self-adjointness-preserving aspect
is contained in Sch\"{u}tzenberger's theorem \cite[Thm. 7.1, p. 18]{BerReu}.
Methodology (iii) has an obscure origin since random walk on many types of graphs has been studied in probability theory for decades and many methods for getting recursions have become commonplace.
In this case we  point to the examples \cite{Aomoto}, \cite{Lalley1},
\cite{NagWoess} and \cite{Woess} as inspirations, and refer the reader to \cite{WoessBook} for background.
\subsubsection{Remark} The form of Proposition \ref{Proposition:MainResultBlockManipulation} is in key respects 
quite similar to 
that of \cite[Thm 2.2]{BelMaiSpei}. Indeed it could not be essentially different since it has the same origins in operator-valued free probability theory.  In particular, \eqref{equation:GSDE2}  can be rewritten as a fixed point equation.
But the matrix upper half-plane plays no role either in the statement or the proof
of Proposition \ref{Proposition:MainResultBlockManipulation}, greatly simplifying matters.

\subsubsection{Remark}
In the semicircular case $\kappa_j^{(\theta)}=\kappa_j(\mu_{x_\theta})=0$ for $j>2$ remarked upon above,
Proposition \ref{Proposition:MainResultAlgebraicityProof} follows straightforwardly
from Propositions \ref{Proposition:AlgebraicityCriterion}.
In  particular, in this relatively simple case,
one can check the nondegeneracy hypothesis of Proposition \ref{Proposition:AlgebraicityCriterion}
directly by using \eqref{equation:GSDE3}.  This observation begins to explain the significance of condition \eqref{equation:GSDE3}.  
A similar approach proves Proposition \ref{Proposition:MainResultAlgebraicityProof} in the more general
case that $\kappa_j^{(\theta)}=\kappa_j(\mu_{x_\theta})=0$ for all but finitely many pairs $(j,\theta)$.

 \subsubsection{Remark}
 We continue in the setup of the preceding remark.
Consider now the remaining case in which $\kappa_j^{(\theta)}=\kappa_j(\mu_{x_\theta})\neq 0$ for infinitely many pairs $(j,\theta)$.
Then one can no longer prove Proposition \ref{Proposition:MainResultAlgebraicityProof} 
by applying Proposition \ref{Proposition:AlgebraicityCriterion} directly to solutions of the system of
equations \eqref{equation:GSDE2} because this system no longer consists exclusively of polynomial equations. This obstruction
is the main difficulty of the proof of Theorem \ref{Theorem:MainResult}. We overcome the obstruction by means of 
the theory of Newton polygons  (see \S\ref{section:NewtonPolygon} below)
and the Weierstrass Preparation Theorem (see \S\ref{section:Gizmo} below). Using these methods
we are able to exhibit a system of $n^2+3qn$ polynomial equations in $n^2+3qn$ unknowns which in a useful sense contains the system  \eqref{equation:GSDE2} and to which Proposition \ref{Proposition:AlgebraicityCriterion} may be applied directly. (See \S\ref{section:Endgame} below.) This enlarged system of equations 
 may be of independent interest. (See \eqref{equation:TheSystem} below.)

\section{Hessenberg-Toeplitz matrices and free cumulants}\label{section:Gentle}
We  introduce the formal version of operator theory used in this paper
and then as an illustration we revisit a  key insight of Voiculescu concerning the free cumulants.

\subsection{The algebras $\Mfrak$ and $\Mfrak((1/z))$}
We introduce two algebras of infinite matrices,
the first an algebra over $\CC$ and the second a larger algebra over $\CC((1/z))$ possessing Banach algebra structure.

 \subsubsection{Notation}
Let $\NN$ denote the set of \textbf{nonnegative} integers. 
\subsubsection{The algebra $\Mfrak$}
Let $\Mfrak$ denote the vector space over $\CC$ consisting of  \linebreak $\NN$-by-$\NN$ matrices $M$ such that for each $j\in \NN$
there exist only finitely many $i\in \NN$ such that $M(i,j)\neq 0$. Note that for each $M\in \Mfrak$
the entry in the upper left corner is denoted by $M(0,0)$, rather than $M(1,1)$, because we are observing
the convention that $0\in \NN$.
Note that every upper-triangular  $\NN$-by-$\NN$ matrix with entries in $\CC$ belongs to $\Mfrak$.
Informally, $\Mfrak$ consists of the ``almost upper-triangular'' matrices.
It is easy to see that matrix multiplication is well-defined on $\Mfrak$ and moreover associative,
thus making $\Mfrak$ into a unital associative algebra with scalar field $\CC$.
Indeed, it is clear that $\Mfrak$ is a copy of the algebra of linear endomorphisms of a complex vector space of countably infinite dimension. We write $\one=1_\Mfrak$ to abbreviate notation.
We equip $\Mfrak$ with the  state $\phi(M)=M(0,0)$, thus defining a noncommutative
probability space $(\Mfrak,\phi)$.  
\subsubsection{The algebra $\Mfrak((1/z))$}
Let $\Mfrak((1/z))$ denote the set of $\NN$-by-$\NN$ matrices $M$ with entries in $\CC((1/z))$
satisfying one and hence both of the following equivalent conditions:
\begin{itemize}
\item There exists a {\em Laurent expansion}
$M=\sum_{n\in \ZZ}M_nz^n$ with coefficients \linebreak $M_n\in \Mfrak$
such that $M_n=0$ for $n\gg 0$.
\item One has $\lim_{i\rightarrow\infty}\val\, M(i,j)=-\infty$
for each $j\in \NN$ (without any requirement of uniformity in $j$) and furthermore one has
$\sup_{i,j\in \NN} \val\, M(i,j)<\infty$.
\end{itemize}
From the equivalent points of view described above it is clear that $\Mfrak((1/z))$ becomes a unital $\CC((1/z))$-algebra
with respect to the usual notion of matrix multiplication.
For $M\in \Mfrak((1/z))$ we define $\val\, M=\sup_{i,j\in\NN} \val\, M(i,j)$.
With respect to the valuation function $\val$ thus extended to $\Mfrak((1/z))$,  the latter becomes a unital
Banach algebra over $\CC((1/z))$. We write $\one=1_\Mfrak=1_{\Mfrak((1/z))}$.

\subsubsection{Elementary matrices and an abuse of notation} Let $\ebold[i,j]\in \Mfrak$ denote the elementary matrix
with entries given by the rule $$\ebold[i,j](k,\ell)=\delta_{ik}\delta_{j\ell}\;\;\mbox{for $i,j,k,\ell\in \NN$.}
$$
The notation $\ebold[i,j]\in \Mfrak$ introduced here is intended to contrast with the notation $\ebold_{ij}\in \Mat_{k\times \ell}(\CC)$
previously introduced for elementary matrices with finitely many rows and columns.
For $M\in \Mfrak$ supported in a set $S\subset \NN\times \NN$ which intersects each column $\NN\times \{j\}$
in a finite set, we abuse
notation by writing 
$$M=\sum_{(i,j)\in S}M(i,j)\ebold[i,j]$$ as a shorthand to indicate the placement of entries of $M$.

 The following simple lemma is a key motivation for the definition of $\Mfrak((1/z))$.
\begin{Lemma}\label{Lemma:ResolventEasy}
Fix $M\in \Mfrak$ arbitrarily and let $\mu$ denote the law of $M$.
Then the matrix 
$z\one-M\in \Mfrak((1/z))$ is invertible
and 
$S_\mu(z)=(z\one-M)^{-1}(0,0)$.
\end{Lemma}
\proof
One has
$$
(z\one-M)^{-1}=\frac{1}{z}\sum_{k=0}^\infty \frac{M^k}{z^k}\in \Mfrak((1/z)).
$$
The geometric series here is convergent because $\val\, \frac{M}{z}<0$. This noted, it is clear that the series $(z\one-M)^{-1}(0,0)$ is the Stieltjes transform of the law of $M$.
\qed

\subsubsection{Remark}
Our setup is inspired by (but is much simpler than) that of \cite{ADKwedge},
and thus belongs to the lineage of \cite{Tate}. The theory of the
$R$-transform overlaps in an interesting way with the theory of
residues developed in \cite{Tate},
one point of contact being the notion of a Hessenberg-Toeplitz matrix. 
(See immediately below.)  This connection deserves further investigation.

\subsection{Hessenberg-Toeplitz matrices}\label{subsection:HTM}
\subsubsection{Basic definitions}
Let $\{\kappa_j\}_{j=1}^\infty$ be any sequence of complex numbers.
Consider the infinite matrix
\begin{equation}\label{equation:HessenbergToeplitzExample}
H=\left[\begin{array}{cccccccccc}
\kappa_1&\kappa_2&\kappa_3&\dots\\
1&\kappa_1&\kappa_2&\kappa_3&\dots\\
&1&\kappa_1&\kappa_2&\kappa_3&\dots\\
&&1&\kappa_1&\kappa_2&\kappa_3&\dots\\
&&&\ddots&\ddots&\ddots&\ddots\end{array}\right]\in \Mfrak.
\end{equation}
Equivalently, in terms of the elementary matrices $\ebold[i,j]\in \Mfrak$ we have
\begin{equation}\label{equation:HessenbergToeplitzExampleBis}
H=\sum_{k\in \NN}\left(\ebold[1+k,k]+\sum_{j\in \NN}\kappa_{j+1}\ebold[k,j+k]\right).
\end{equation}
The matrix $H$ displays the (upper) {\em Hessenberg} pattern: $i>j+1\Rightarrow H(i,j)=0$ for $i,j\in \NN$.
The matrix $H$ also displays the {\em Toeplitz} pattern: 
$H(i+1,j+1)=H(i,j)$ for $i,j\in \NN$.
Accordingly we call $H$ a {\em Hessenberg-Toeplitz matrix}.

The reason for our interest in the matrix $H$ is explained by the next lemma.
\begin{Lemma}\label{Lemma:HessenbergToeplitz}
Assumptions and notation are as above. Then for every positive integer $j$,
the $j^{th}$ free cumulant of $H$ viewed as a noncommutative random variable in the noncommutative
probability space $(\Mfrak,\phi)$ equals $\kappa_j$.
\end{Lemma}
 This fact is well-known---it is a key insight for Voiculescu's theory of the \linebreak $R$-transform 
\cite{Voiculescu}.
It is therefore not  necessary to give a proof. But we nevertheless give a proof in \S\ref{subsection:ProofOfHT} below  
by manipulation of block-decomposed matrices
in order to set the pattern for the more difficult calculations undertaken in \S\ref{section:GSDE}  below.
\subsubsection{Remark} If one patiently works through the definitions
and uses the theory of residues from \cite{Tate}, one can see that
Lemma \ref{Lemma:HessenbergToeplitz} is equivalent 
to the one-variable Lagrange Inversion Formula. 

\subsection{Inversion of block-decomposed matrices}
We pause to review a method of calculation used repeatedly in the sequel.
\subsubsection{Basic identities}
Let
$$\left[\begin{array}{cc}
\abold&\bbold\\
\cbold&\dbold\end{array}\right]$$
be an invertible square matrix (in practice infinite)
decomposed into blocks where $\abold$ and $\dbold$ are square
and $\dbold$ is also invertible. Then we have a factorization
\begin{equation}\label{equation:PreBigInversion}
\left[\begin{array}{cc}
\abold-\bbold \dbold^{-1}\cbold&\zero\\
\zero&\dbold
\end{array}\right]
=\left[\begin{array}{cc}
\one&-\bbold \dbold^{-1}\\
\zero&\one\end{array}\right]\left[\begin{array}{cc}
\abold&\bbold\\
\cbold&\dbold
\end{array}\right]\left[\begin{array}{cc}
\one&\zero\\
-\dbold^{-1}\cbold&\one
\end{array}\right]
\end{equation}
from which in particular we infer that the {\em Schur complement} 
$\abold-\bbold \dbold^{-1}\cbold$ is invertible.
Let 
$$\gbold=(\abold-\bbold\dbold^{-1}\cbold)^{-1}.$$
From \eqref{equation:PreBigInversion} one then straightforwardly derives
the inversion formula
\begin{eqnarray}\label{equation:BigInversion}
\left[\begin{array}{cc}
\abold&\bbold\\
\cbold&\dbold\end{array}\right]^{-1}&=&
\left[\begin{array}{cc}
\zero&\zero\\
\zero&\dbold^{-1}
\end{array}\right]+\left[\begin{array}{r}
\one\\
-\dbold^{-1}\cbold
\end{array}\right]
\gbold\left[\begin{array}{rr}
\one&-\bbold \dbold^{-1}
\end{array}\right].
\end{eqnarray}
The latter formula also shows that invertibility of $\dbold$ and $\abold-\bbold\dbold^{-1}\cbold$
implies invertibility of $\left[\begin{array}{cc}
\abold&\bbold\\
\cbold&\dbold\end{array}\right]$.
For convenient  application in \S\ref{section:GSDE}, we restate in abstract form a couple of relations
among blocks recorded in formula \eqref{equation:BigInversion}.

\begin{Lemma}\label{Lemma:AbstractSchurComplement}
Let $\AAA$ be a unital associative algebra (perhaps not commutative).
Let $\pi,\sigma\in \AAA$ satisfy $\pi^2=\pi\neq 0$, $\sigma^2=\sigma\neq 0$, $\pi\sigma=\sigma \pi=0$
and $1_\AAA=\pi+\sigma$. Let $A\in \AAA$ be invertible.
Assume furthermore that $\sigma A \sigma$ is invertible in the algebra $\sigma \AAA\sigma$
and let $A_\sigma^{-1}$ denote the inverse.
Then we have
\begin{eqnarray}
\label{equation:SC3}
\sigma A^{-1}\pi&=&-A_\sigma^{-1} A\pi A^{-1}\pi\;\;\mbox{and}\\
\label{equation:SC4}
A^{-1}\sigma&=&(1_\AAA-A^{-1}\pi A\sigma)A_\sigma^{-1}.
\end{eqnarray}
\end{Lemma}
\proof  We have
$$\sigma A\sigma A^{-1}\pi=-\sigma AA^{-1}\pi+\sigma A\sigma A^{-1}\pi=-\sigma A\pi A^{-1}\pi.$$
Now left-multiply extreme terms by $A_\sigma^{-1}$ to recover \eqref{equation:SC3}. Similarly, we have
$$\sigma=\sigma AA_\sigma^{-1}=(A -\pi A \sigma)A_\sigma^{-1}.$$
Now left-multiply extreme terms by $A^{-1}$ to recover \eqref{equation:SC4}. \qed

\subsection{Proof of Lemma \ref{Lemma:HessenbergToeplitz}}\label{subsection:ProofOfHT}
 Consider the Laurent series
$$f=f(z)=z+\sum_{j=1}^\infty \kappa_jz^{j-1}\in z+\CC[[1/z]].$$
It will suffice to show that $f$ is equal to the modified $R$-transform of the law of $H$.
Consider also
the Stieltjes transform
$$g=g(z)=S_{\mu_H}(z)\in (1/z)+(1/z^2)\CC[[1/z]]$$
of the law of $H$. Since $z+\CC[[1/z]]$ forms a group under composition, 
it will suffice to show that 
$z=f\circ\frac{1}{g}$,
equivalently $z=g^{-1}+\sum_{j=1}^\infty \kappa_j g^{j-1}$,
or equivalently
$$
1=(z-\kappa_1)g-\sum_{j=2}^\infty \kappa_j g^{j}.
$$
Let 
$$A=z\one-H\in \Mfrak((1/z)).$$
By Lemma \ref{Lemma:ResolventEasy}
the inverse 
$$G=A^{-1}\in \Mfrak((1/z))$$
exists
and furthermore
$$g=G(0,0).$$
In view of the relation 
$$1=\sum_{k\in \NN}A(0,k)G(k,0)$$
holding because $G=A^{-1}$,
it will be enough simply to prove that
\begin{equation}\label{equation:HessTopClincherBis}
G(i,0)=g^{i+1}\;\;\mbox{for $i\in \NN$.}
\end{equation}
Now with an eye toward applying \eqref{equation:BigInversion} above,
  consider the block decomposition
$$A=\left[\begin{array}{cc}
\abold&\bbold\\
\cbold&\dbold\end{array}\right]$$
where
$$\abold=z-\kappa_1,\;\;\bbold=-\left[\begin{array}{cccc}
\kappa_2&\kappa_3&\dots\end{array}\right],\;\;
\cbold=-\left[\begin{array}{c}1\\
0\\
\vdots\end{array}\right],\;\;\mbox{and}\;\;\dbold=A.$$
By \eqref{equation:BigInversion} we have
$$
\left[\begin{array}{c}G(1,0)\\
G(2,0)\\
\vdots
\end{array}\right]=-\dbold^{-1}\cbold\gbold=
\left[\begin{array}{cccccc}
G(0,0)\\G(1,0)\\ \vdots
\end{array}\right]g,$$
whence \eqref{equation:HessTopClincherBis}. The proof is complete.
 \qed
 \section{The linearization step}
\label{section:Linearization}
In this section we apply (a fragment of) Sch\"{u}tzenberger's Theorem to a suitable model of free noncommutative random variables
with prescribed free cumulants, thus advancing the proof of Proposition \ref{Proposition:MainResultBlockManipulation}.
\subsection{Stars and diamonds}
We build a model for the free unital associative monoid on $q$ generators for which
$\NN$ is the underlying set. Using this monoid structure we will be able to construct and manipulate usefully patterned matrices in $\Mfrak((1/z))$.
\subsubsection{Notation}
The parameter $q$ figuring in Theorem \ref{Theorem:MainResult}
is considered fixed throughout the remainder of the paper. Many constructs below depend on $q$ but the notation frequently does not show this.
 
\subsubsection{Improper representations to the base $q$}
Suppose at first that $q>1$.
In grade school one learns to represent nonnegative integers to the base $q$
using place-notation and digits selected from the set $\{0,\dots,q-1\}$.
It is not hard to see that using instead digits selected from the set $\{1,\dots,q\}$
one still gets a unique representation for every member of $\NN$,
it being understood that $0$ is represented by the empty digit string $\emptyset$.
A representation to the base $q$ of a nonnegative integer using digits $\{1,\dots,q\}$ will be called {\em improper}.
Improper representations to the base $q$ make sense also for $q=1$. In the latter extreme case
each $x\in \NN$ is represented by a string of $1$'s of length $x$.

\subsubsection{Example: counting improperly to the base $3$}
$$\emptyset, 1,2,3,11,12,13,21,22,23,31,32,33,111,112,113,121,122,\dots$$

\subsubsection{The binary operation $\star_q$}
We define the binary operation
$$\star=\star_q:\NN\times \NN\rightarrow\NN$$
by the rule 
$$x\star y=xq^\ell+y\;\;\mbox{where $\ell\in\NN$ satisfies}\;\; \frac{q^\ell-1}{q-1}\leq y\leq q\cdot\left(\frac{q^\ell-1}{q-1}\right).$$
Informally, $\ell$ is the number of digits in the improper base $q$ representation of $y$,
and  $x\star y$ is the concatenation of the digit strings of $x$ and of $y$.
The operation $\star$ makes $\NN$ into a free associative monoid freely generated by the digits $1,\dots,q$
with $0$ as the identity element.
\begin{Lemma}\label{Lemma:Unpin}
$\NN\setminus \{0\}$ is the disjoint union of the sets $\NN\star \theta$ for $\theta=1,\dots,q$.
\end{Lemma}
 There is nothing to prove. We record this for convenient reference since however trivial,
this is an important point in a later proof.
\subsubsection{The binary operation $\diamond_q$}
We define the binary operation 
$$\diamond=\diamond_q:\NN\times \NN\rightarrow \NN$$ by the formula
$$x^{\diamond y}=\underbrace{x\star\cdots \star x}_y.$$
 We use  exponential-style notation
to emphasize the analogy with exponentiation in the usual sense. 
\begin{Lemma}\label{Lemma:UniqueFactorization}
For $\theta=1,\dots,q$, every $x\in \NN$ has a unique factorization $x=\theta^{\diamond i}\star k$
where $i\in \NN$ and $k\in \NN\setminus \theta\star\NN$.
\end{Lemma}
 Again, there is nothing to prove. We record this  for convenient reference.

\subsubsection{Remark}\label{subsubsection:TreePicture}
Consider the  graph $\Gamma=\Gamma_q$ with vertex set $\NN$
and an edge connecting $x$ to $\theta\star x$ for each
pair $(\theta,x)\in \{1,\dots,q\}\times \NN$.
With $0\in \NN$ designated as the root, the resulting graph $\Gamma$ 
is an infinite rooted planar tree in which every vertex has a ``birth-ordered'' set of $q$ children, i.e., a $q$-ary rooted tree.
We do not explicitly use the $q$-ary tree in this paper
because we instead rely on the monoid $(\NN,\star)$ to do our bookkeeping.
Nonetheless the notion of the $q$-ary tree remains a strong guide to intuition.

\subsection{Free random variables with prescribed free cumulants}
The next proposition exhibits a model for $q$ free noncommutative random variables with prescribed free cumulants.
The model is essentially the same as that used in \cite{Voiculescu}, but with the notation designed to make recursions easy to see.

\begin{Proposition}\label{Proposition:EmbeddedHT}
Let 
$$\left\{\left\{\kappa_j^{(\theta)}\right\}_{j=1}^\infty\right\}_{\theta=1}^q$$ be any family of complex numbers.
Then the family
\begin{equation}\label{equation:ExcellentFamily}
\sum_{k\in \NN}\ebold [\theta\star k,k]+\sum_{j\in \NN}\sum_{k\in \NN}\kappa_{j+1}^{(\theta)}\ebold[ k,\theta^{\diamond j}\star k]\in\Mfrak\;\;\mbox{for $\theta=1,\dots,q$}
\end{equation}
of noncommutative random variables is freely independent
and moreover the $j^{th}$ free cumulant of the $\theta^{th}$ noncommutative random variable equals $\kappa_j^{(\theta)}$.
\end{Proposition}
\noindent The proof requires some preparation and is completed in \S\ref{subsubsection:ModelProofEnd}.

\subsubsection{Self-embeddings of $\Mfrak$}

For $\theta=1,\dots q$ and $A\in \Mfrak$ we define
\begin{equation}\label{equation:AutoEmbedding}
A^{(\theta)}=\sum_{k\in \NN\setminus \theta\star \NN}\;\;
\sum_{i,j\in \NN} A(i,j)\ebold[\theta^{\diamond i}\star k,\theta^{\diamond j}\star k].\end{equation}
By Lemma \ref{Lemma:UniqueFactorization} the matrix $A^{(\theta)}$ is block-diagonal with copies of $A$ indexed by $\NN\setminus \theta\star \NN$
repeated along the diagonal.
Thus the map $(A\mapsto A^{(\theta)}):\Mfrak\rightarrow\Mfrak$
is a unital one-to-one homomorphism of algebras. Note that $A^{(\theta)}(0,0)=A(0,0)$
and hence the map $A\mapsto A^{(\theta)}$ is law-preserving. Let $\Mfrak^{(\theta)}$
denote the embedded image of $\Mfrak$ under the map $A\mapsto A^{(\theta)}$.

\begin{Lemma}\label{Lemma:SelfFreeing}
The subalgebras $\Mfrak^{(1)},\dots,\Mfrak^{(q)}\subset \Mfrak$ are freely independent.
\end{Lemma}
\proof Fix $\theta_1,\dots,\theta_k\in \{1,\dots,q\}$ such that 
$$\theta_1\neq \theta_2,\;\; \theta_2\neq \theta_3\;\;, \dots, \;\;\theta_{k-1}\neq \theta_k.$$
Fix $A_1,\dots,A_k\in \Mfrak$ such that 
$$A_1(0,0)=\cdots=A_k(0,0)=0.$$ Our task is to verify that
\begin{eqnarray}\label{equation:FreeIndependenceNuff}
&&\sum_{(i_1,\dots,i_{k-1})\in \NN^{k-1}}A_1^{(\theta_1)}(0,i_1)A_2^{(\theta_2)}(i_1,i_2)
\cdots A_{k-1}^{(\theta_{k-1})}(i_{k-2},i_{k-1})A_k^{(\theta_k)}(i_{k-1},0)\\
\nonumber&=&(A_1^{(\theta_1)}\cdots A_k^{(\theta_k)})(0,0)=0.
\end{eqnarray}
Now by definition, for any matrix $A\in \Mfrak$ such that $A(0,0)=0$ and $\theta=1,\dots,q$, 
the matrix entry $A^{(\theta)}(i,j)$ vanishes unless $i\in \theta\star\NN$ or $j\in \theta\star \NN$.
It follows that for any \linebreak $(k-1)$-tuple $(i_1,\dots,i_{k-1})\in \NN^{k-1}$ the corresponding term
in the sum on the left side of \eqref{equation:FreeIndependenceNuff} vanishes.
Thus, {\em a fortiori},  \eqref{equation:FreeIndependenceNuff} holds.
\qed
\subsubsection{Proof of Proposition \ref{Proposition:EmbeddedHT}}\label{subsubsection:ModelProofEnd}
Let $H_\theta$ be a copy of the matrix $H$ defined in \eqref{equation:HessenbergToeplitzExample}
and equivalently in \eqref{equation:HessenbergToeplitzExampleBis}, with $\kappa_j$ replaced by $\kappa_j^{(\theta)}$.
By Lemma \ref{Lemma:HessenbergToeplitz} we know already that the $j^{th}$ free cumulant of the law of $H_\theta$ equals $\kappa_j^{(\theta)}$.
Substituting directly into the definition \eqref{equation:AutoEmbedding} we  have
$$
H^{(\theta)}_\theta= \sum_{i\in\NN}\sum_{k\in \NN\setminus \theta\star \NN}\ebold[\theta^{\diamond (i+1)}\star k,\theta^{\diamond i}\star k]
+ \sum_{i\in \NN}\sum_{j\in \NN}  \sum_{k\in \NN\setminus \theta\star \NN}\kappa_{j+1}^{(\theta)}\ebold[\theta^{\diamond i}\star k,\theta^{\diamond (i+j)}\star k].
$$
The result follows now via
Lemmas \ref{Lemma:UniqueFactorization} and \ref{Lemma:SelfFreeing}. \qed

\subsubsection{Remark} \label{Remark:LoweringAndRaising}
Voiculescu \cite{Voiculescu} introduced the Boltzmann-Fock space model  of free random
variables using lowering and raising operators  for his striking proof of additivity of the $R$-transform 
for addition of free random variables.  Also see \linebreak \cite[Cor. 5.3.23]{AGZ} and its proof
for a quick review of this material.
Proposition \ref{Proposition:EmbeddedHT} is merely a description
of the Boltzmann-Fock space model using notation chosen to make recursions more easily accessible.
In the setup of Proposition \ref{Proposition:EmbeddedHT} the matrices
\begin{equation}\label{equation:Sinister}
\hat{\lambda}^{(\theta)}=\sum_{k\in \NN}\ebold[\theta\star k,k]\in \Mfrak\;\;\mbox{and}\;\;
\lambda^{(\theta)}=\sum_{k\in \NN}\ebold[k,\theta\star k]\in \Mfrak\;\;
\mbox{for $\theta=1,\dots,q$}
\end{equation}
correspond to the lowering and raising operators considered in \cite{Voiculescu}, respectively.
Note that using the operators on line \eqref{equation:Sinister}
 we can rewrite the operators on \eqref{equation:ExcellentFamily} in the more familiar form
\begin{equation}\label{equation:ExcellentFamilyBis}
\hat{\lambda}^{(\theta)}+\sum_{j\in \NN}\kappa_{j+1}^{(\theta)}(\lambda^{(\theta)})^j\in\Mfrak\;\;\mbox{for $\theta=1,\dots,q$.}
\end{equation}
The infinite sum here is an abuse of notation but it nonetheless makes sense because the matrices being summed have disjoint supports only finitely many of which meet any given column.
Later  we will also consider lowering and raising operators
\begin{equation}\label{equation:Dextrous}
\hat{\rho}^{(\theta)}=\sum_{k\in \NN}\ebold[k\star \theta,k]\in \Mfrak\;\;\mbox{and}\;\;
\rho^{(\theta)}=\sum_{k\in \NN}\ebold[k, k\star \theta]\in \Mfrak\;\;
\mbox{for $\theta=1,\dots,q$}
\end{equation}
acting (so to speak) on the right rather than the left. 
The relations
\begin{eqnarray}\label{equation:Commutation1}
&&\lambda^{(\theta)}\rho^{(\theta')}=\rho^{(\theta')}\lambda^{(\theta)},\;\;
\hat{\lambda}^{(\theta)}\hat{\rho}^{(\theta')}=
\hat{\rho}^{(\theta')}\hat{\lambda}^{(\theta)},\\
\label{equation:Commutation2}
&&\rho^{(\theta')}\hat{\rho}^{(\theta)}=\delta_{\theta\theta'}\one,\;\;\mbox{and}\;\;
\sum_{\alpha=1}^q\hat{\rho}^{(\alpha)}\rho^{(\alpha)}=\sum_{i=1}^\infty \ebold[i,i]
\end{eqnarray}
for $\theta,\theta'=1,\dots,q$ are easy to verify. 
See for example  \cite[Sec. 3.4]{Anderson} where these and further relations are written out as part of an analysis
leading (without any reference to noncrossing partitions) to the Schwinger-Dyson equation for semicircular variables.
\subsubsection{Remark}
The interplay of left and right lowering and raising operators is a fundamental feature of the recently introduced bi-free framework of \cite{VoiculescuJanus}.

\subsection{Kronecker products and the isomorphism $\natural$}
We  introduce notation which is rather tedious to define but convenient to calculate with.
 \subsubsection{Classical Kronecker products}
 Recall that for matrices of finite size the {\em Kronecker product} 
$$A^{(1)}\otimes A^{(2)}\in \Mat_{k_1k_2\times \ell_1\ell_2}(\CC)\;\;\;
\left(A^{(\alpha)}\in \Mat_{k_\alpha\times \ell_\alpha}(\CC)\;\;\mbox{for $\alpha=1,2$}\right)$$  is defined by the rule
\begin{equation*}
A^{(1)}\otimes A^{(2)}=\left[\begin{array}{ccc}
A^{(1)}(1,1)A^{(2)}&\dots&A^{(1)}(1,\ell_1)A^{(2)}\\
\vdots&&\vdots\\
A^{(1)}(k_1,1)A^{(2)}&\dots&A^{(1)}(k_1,\ell_1)A^{(2)}
\end{array}\right]
\end{equation*}
or equivalently and more explicitly (if more cumbersomely)
\begin{eqnarray*}
&&A^{(1)}\otimes A^{(2)}(n(i_1-1)+i_2,n(j_1-1)+j_2)=A^{(1)}(i_1,j_1)A^{(2)}(i_2,j_2)\\
&&\hspace{1.8in}\;\mbox{for $\alpha=1,2$, $i_\alpha=1,\dots,k_\alpha$, and $j_\alpha=1,\dots,\ell_\alpha$.}
\end{eqnarray*}
\subsubsection{Kronecker products involving infinite matrices}
In the mixed infinite/finite case we define the {\em Kronecker product}
$$x\otimes a\in \Mfrak((1/z))\;\;\;\left(x\in \Mfrak((1/z))\;\;\mbox{and}\;\;a\in \Mat_n(\CC((1/z)))\right)$$
by the rule
$$x\otimes a=\left[\begin{array}{ccccc}
x(0,0)a&x(0,1)a&\dots\\
x(1,0)a&x(1,1)a&\dots\\
\vdots&\vdots&\ddots
\end{array}\right]$$
or equivalently and more explicitly
\begin{eqnarray*}
&&(x\otimes a)(i_1n+i_2-1,j_1n+j_2-1)=x(i_1,j_1)a(i_2,j_2)\\
\nonumber&&\hspace{2.0in}\mbox{for $i_1,j_1\in \NN$ and $i_2,j_2=1,\dots,n$.}
\end{eqnarray*}
We also define 
$$a\otimes x\in \Mat_n(\Mfrak((1/z)))\;\;\;\left(a\in \Mat_n(\CC((1/z)))\;\;\mbox{and}\;\;
x\in \Mfrak((1/z))\right)$$
 by the somewhat ungainly iterated index formula
 \begin{eqnarray*}
&&((a\otimes x)(i_2,j_2))(i_1,j_1)=a(i_2,j_2)x(i_1,j_1)\\
\nonumber&&\hspace{2.0in}\mbox{for $i_1,j_1\in \NN$ and $i_2,j_2=1,\dots,n$.}
\end{eqnarray*}

\subsubsection{The operation $\natural$}
For $M\in \Mat_n(\Mfrak((1/z)))$ we define $M^\natural\in \Mfrak((1/z))$
by the formula
$$M^\natural=\sum_{i_1,j_1\in\NN}\sum_{i_2,j_2=1}^n
(M(i_2,j_2)(i_1,j_1))(\ebold[i_1,j_1]\otimes \ebold_{i_2j_2}),$$
thus defining an isometric isomorphism
$$(M\mapsto M^\natural):\Mat_n(\Mfrak((1/z)))\rightarrow \Mfrak((1/z))$$
of Banach algebras over $\CC((1/z))$, where the source algebra is given
Banach algebra structure by declaring that
$\val\, A=\max_{i,j=1}^n \val\, A(i,j)$
for $A\in \Mat_n(\Mfrak((1/z)))$.   Finally, note that
\begin{equation}\label{equation:Switcharoo}
(a\otimes x)^\natural=x\otimes a\;\;\;\mbox{for $a\in \Mat_n(\CC((1/z)))$ and $x\in \Mfrak((1/z))$}.
\end{equation}
Thus the operation $\natural$ has a natural interpretation as exchange of tensor factors.
\begin{Lemma}\label{Lemma:FormalResolvent}
Fix $A\in \Mat_n(\Mfrak)$.
Then the following statements hold:
\begin{eqnarray*}
&&\mbox{$zI_n\otimes \one-A\in \Mat_n(\Mfrak((1/z)))$ and $z\one-A^\natural\in \Mfrak((1/z))$ are invertible.}\\
&&\left((zI_n\otimes \one-A)^{-1}\right)^\natural=(z\one-A^\natural)^{-1}.\\
&&S_{\mu_A}(z)=\frac{1}{n}\sum_{i=0}^{n-1}(z-A^\natural)^{-1}(i,i)
=\left(\frac{1}{n}\sum_{i=1}^n(zI_n\otimes \one-A)^{-1}(i,i)\right)(0,0).
\end{eqnarray*}
\end{Lemma}
 This statement supplements Lemma \ref{Lemma:ResolventEasy} only by some minor bookkeeping details.
We therefore omit proof.

\subsection{Digital linearization}
Here is the main result in this section.
\begin{Proposition}\label{Proposition:DigitalLinearization}
Let $(\AAA,\phi)$ be a noncommutative probability space.
Let \linebreak
$x_1,\dots,x_q\in \AAA$ be freely independent noncommutative random variables.
Fix \linebreak
$X\in\Mat_p(\CC\langle x_1,\dots,x_q\rangle)\subset\Mat_p(\AAA)$.
Let  $\kappa_j^{(\theta)}=\kappa_j(\mu_{x_{\theta}})$ for $j=1,2,\dots$ and $\theta=1,\dots,q$.
Then there exist for some $N>0$ matrices $L_0,L_1,\dots,L_q\in \Mat_{p+N}(\CC)$
with the following properties:
\begin{eqnarray}\label{equation:L0support}
&&\mbox{$L_0$ vanishes identically in the upper left $p$-by-$p$ block.}\\
\label{equation:L+support}
&&\mbox{$L_1,\dots,L_q$ are supported in the lower right $N$-by-$N$ block.}\\
\label{equation:Ldef}
&&L=\one\otimes \left(L_0+\left[\begin{array}{cc}
zI_p&0\\
0&0\end{array}\right]\right)+\sum_{\theta=1}^q\sum_{k\in \NN}\ebold [\theta\star k,k]\otimes L_\theta\\
\nonumber&&+\sum_{\theta=1}^q\sum_{j\in \NN}\sum_{k\in \NN}\kappa^{(\theta)}_{j+1}\ebold[ k,\theta^{\diamond j}\star k]\otimes L_\theta\in\Mfrak((1/z))\;\;\mbox{is invertible.}\\
\label{equation:GeneralizedResolventKey}
&&S_{\mu_X}(z)=\frac{1}{p}\sum_{i=0}^{p-1} L^{-1}(i,i).
\end{eqnarray}
\end{Proposition}
 We call $L$ a {\em digital linearization} of $X$. It is worth remarking that this linearization
is thoroughgoing in the sense that not only do the variables $x_1,\dots,x_q$ appear linearly---so also does the variable $z$.
The proof will be completed in \S\ref{subsection:ProofOfDigitalLinearization} below.

\subsubsection{Remark}\label{subsubsection:TreePictureBis}
Picking up again on the idea mentioned in \S\ref{subsubsection:TreePicture},
and adopting the absurd point of view that probabilities can be square matrices with arbitrary complex number entries,
the matrix $L$ describes a random walk
on the $q$-ary tree $\Gamma_q$ such that from a given vertex $x\in \NN$,
one may  (i) step one unit back toward the root (if not already at the root), (ii) stay in place,
or (iii) step away from the root arbitrarily far along along a geodesic $\{\theta^{\diamond i}\star x\mid i\in \NN\}$
for some \linebreak $\theta\in \{1,\dots,q\}$. Whether or not this interpretation of $L$ is absurd, it is does make
random walk intuition available to analyze $L$. Guided by this intuition we will prove in \S\ref{section:GSDE} below
that the generalized Schwinger-Dyson equation holds for ``upper left corners'' of  matrices of the form $L$,
as well as for more general infinite matrices.

\subsection{Sch\"{u}tzenberger's theorem}
The next lemma recalls what we need of the self-adjoint linearization trick,
and as we have already noted in the introduction, what we need boils down to Sch\"{u}tzenberger's Theorem \cite[Thm. 7.1]{BerReu}. In fact we need only a quite specialized consequence of this theorem, 
or rather of its proof, simple enough to prove quickly from scratch, as follows. 
\begin{Lemma}\label{Lemma:Linearization}
 For each $f\in\Mat_p(\CC\langle \Xbold_1,\dots,\Xbold_q\rangle)$ there exists a factorization $f=bd^{-1}c$ (called a {\em linearization} of $f$) where 
$$b\in \Mat_{p\times N}(\CC),\;\;
c\in \Mat_{N\times p}(\CC), \;\;d\in\GL_N(\CC\langle \Xbold_1,\dots,\Xbold_q\rangle),$$
and each entry of $d$ belongs to the $\CC$-linear span of $1,\Xbold_1,\dots,\Xbold_q$.
 \end{Lemma}
 \noindent Note that the proof below actually produces $d$ with the further property that $d-I_N$ is strictly upper triangular.
 \proof If every entry of $f$ belongs to the $\CC$-linear
span of $1,\Xbold_1,\dots,\Xbold_q$, then, say, 
$$f=
\left[\begin{array}{cc}
I_p&0\end{array}\right]\left[\begin{array}{cc}
I_p&-f\\
0&I_p\\
\end{array}\right]^{-1}\left[\begin{array}{c}0\\I_p\end{array}\right]$$ is a linearization.
Thus it will be enough to demonstrate that given linearizable $f_1,f_2\in \Mat_p(\CC\langle \Xbold_1,\dots,\Xbold_q\rangle)$,
again $f_1+f_2$ and $f_1f_2$ are linearizable.
So suppose that $f_i=b_id_i^{-1}c_i$ for $i=1,2$ are factorizations of the desired form.
We then have
$$f_1+f_2=\left[\begin{array}{cc}
b_1&b_2\end{array}\right]\left[\begin{array}{cc}
d_1&0\\
0&d_2
\end{array}\right]^{-1}\left[\begin{array}{c}
c_1\\c_2
\end{array}\right]\;\;\mbox{and}\;\;$$
$$f_1f_2=\left[\begin{array}{cccc}
b_1&0&0\end{array}\right]
\left[\begin{array}{cccc}
d_1&c_1&0\\
0&1&b_2\\
0&0&d_2
\end{array}\right]^{-1}\left[\begin{array}{c}
0\\0\\c_2
\end{array}\right].$$
To assist the reader in checking the second formula, we note that
$$
\left[\begin{array}{cccc}
d_1&c_1&0\\
0&1&b_2\\
0&0&d_2
\end{array}\right]^{-1}=\left[\begin{array}{ccc}
d_1^{-1}&-d_1^{-1}c_1&d^{-1}c_1b_2d_2^{-1}\\
0&1&-b_2d_2^{-1}\\
0&0&d_2^{-1}
\end{array}\right].$$
Thus $f_1+f_2$ and $f_1f_2$ have linearizations. Consequently the lemma does indeed hold.
 \qed
\subsubsection{Remark}\label{subsubsection:Automaton}
We refer the reader to the book \cite{BerReu}
for a complete discussion of Sch\"{u}tzenberger's Theorem and its context in the theory of rational formal noncommutative power series. Nonetheless, we feel that we do owe the reader at least a brief sketch of the interpretation of Lemma \ref{Lemma:Linearization}
that identifies it as a consequence of Sch\"{u}tzenberger's Theorem.
For simplicity and with some loss of generality we assume that $p=1$. Without further loss of generality we may assume that
$$d(0)=I_N-\sum_{\theta=1}^qh_\theta\otimes \Xbold_\theta.$$
Let
$$h=(h_1,\dots,h_q)\in \Mat_{N}^q.$$
Then we have
$$d^{-1}= \sum_{\begin{subarray}{c}
\mbox{\scriptsize monomials}\\
\Mbold\in \CC\langle \Xbold\rangle
\end{subarray}}\Mbold(h)\otimes \Mbold
$$
(the sum is actually finite on account of the remark immediately following the statement of 
Lemma \ref{Lemma:Linearization}) and hence
$$f=bd^{-1}c=\sum_{\begin{subarray}{c}
\mbox{\scriptsize monomials}\\
\Mbold\in \CC\langle \Xbold\rangle
\end{subarray}}(b\Mbold(h)c)\Mbold.
$$
Write
$$f=\sum_{\begin{subarray}{c}
\mbox{\scriptsize monomials}\\
\Mbold\in \CC\langle \Xbold\rangle
\end{subarray}} a_\Mbold \Mbold\;\;(a_\Mbold\in \CC).$$
In the spirit of Sch\"{u}tzenberger's theorem 
we should think of the collection
$$(b,h_1,\dots,h_q,c)$$ as 
a ``linear automaton'' which 
by the rule $$a_\Mbold=b\Mbold(h)c\;\;\mbox{for monomials $\Mbold\in \CC\langle \Mbold\rangle$}$$ 
``computes'' $f$ coefficient-by-coefficient. More generally and analogously, \linebreak Sch\"{u}tzenberger's theorem produces a
linear automaton  to compute the coefficients of any given rational noncommutative formal power series,
and moreover asserts that any noncommutative formal power series so ``computable'' is rational.
We note finally that our application of Lemma \ref{Lemma:Linearization} below is not dependent
on its ``automatic'' interpretation. Neither the theory of automata nor the theory of rational formal noncommutative power series are needed in the sequel.

\subsection{Proof of Proposition \ref{Proposition:DigitalLinearization}}\label{subsection:ProofOfDigitalLinearization}
Without loss of generality we may assume that $(\AAA,\phi)$ is the noncommutative probability
space $(\Mfrak,\phi)$ and we may take $\{x_\theta\}_{\theta=1}^q$ to be the family constructed in Proposition \ref{Proposition:EmbeddedHT}.
 Fix
$$f=f(\Xbold_1,\dots,\Xbold_q)\in \Mat_p(\CC\langle \Xbold_1,\dots,\Xbold_q\rangle)$$ such that
$f(x_1,\dots,x_q)=X$,
write  $f=bd^{-1}c$ as in Lemma \ref{Lemma:Linearization},
and then write
\begin{eqnarray*}
\left[\begin{array}{cc}
0&b\\
c&d\end{array}\right]&=&L_0\otimes 1_{\CC\langle \Xbold_1,\dots,\Xbold_q\rangle}+L_1\otimes \Xbold_1+
\cdots+L_q\otimes \Xbold_q\\
&\in &\Mat_{p+N}(\CC\langle \Xbold_1,\dots,\Xbold_q\rangle)\;\;\;(L_0,L_1,\dots,L_q\in \Mat_{p+N}(\CC))
\end{eqnarray*}
in the unique possible way. 
Finally, let
$$L=\one\otimes \left(L_0+\left[\begin{array}{cc}
zI_p&0\\
0&0\end{array}\right]\right)+x_1\otimes L_1+
\cdots+x_q\otimes L_q\in\Mfrak((1/z))$$
noting that this expression when expanded in terms of elementary matrices
takes by \eqref{equation:ExcellentFamily} the desired form \eqref{equation:Ldef}.
Let 
$$B=b\otimes \one\in \Mat_{p\times N}(\Mfrak)\;\;\mbox{and}\;\;C=c\otimes \one\in \Mat_{N\times p}(\Mfrak).
$$
Let
$$D\in \GL_N(\Mfrak)\;\;\mbox{and}\;\;F\in \Mat_p(\Mfrak)$$
be the evaluations of $d$ and $f$, respectively, at $\Xbold_\theta=x_\theta$ for $\theta=1,\dots,q$.
Since $D$ is the image of an invertible matrix under a unital algebra homomorphism, $D$ is indeed invertible.
Furthermore $zI_p\otimes \one-F$ is invertible by Lemma \ref{Lemma:FormalResolvent}.
It follows by the discussion after formula \eqref{equation:BigInversion}
that the matrix
$$\left[\begin{array}{cc}
zI_p\otimes \one&B\\C&D\end{array}\right]\in \Mat_{p+N}(\Mfrak((1/z)))
$$
is invertible,
and from \eqref{equation:BigInversion} itself it follows that
\begin{eqnarray*}
&&\left[\begin{array}{cc}
zI_p\otimes \one&B\\C&D\end{array}\right]^{-1}\\
\nonumber&=&\left[\begin{array}{cc}
0&0\\
0&D^{-1}
\end{array}\right]+\left[\begin{array}{r}
I_p\otimes \one\\
-D^{-1}C
\end{array}\right](zI_p\otimes \one-F)^{-1}\left[\begin{array}{rr}
I_p\otimes \one&-BD^{-1} 
\end{array}\right].
\end{eqnarray*}
In turn, by \eqref{equation:Switcharoo} we have
$$L=\left[\begin{array}{rl}
zI_p\otimes \one&B\\
C&D
\end{array}\right]^\natural,$$
hence $L$ is invertible and moreover \eqref{equation:GeneralizedResolventKey} holds
by Lemma \ref{Lemma:FormalResolvent}.
The proof of Proposition \ref{Proposition:DigitalLinearization} is complete.
 \qed

\section{Solving the generalized Schwinger-Dyson equation}\label{section:GSDE}
We finish the proof of Proposition \ref{Proposition:MainResultBlockManipulation}
by constructing sufficiently many solutions of the generalized Schwinger-Dyson equation.

\subsection{Statement of the construction}
Here is our main result in this section.
\begin{Proposition}\label{Proposition:GSDE}
Fix data of the form \eqref{equation:Data}.
Consider the matrix
\begin{equation}\label{equation:Adef}
A=-\one\otimes a^{(0)}-\sum_{\theta=1}^q\sum_{k\in \NN}
\left(\ebold[\theta\star k,k]+\sum_{j=1}^\infty\kappa_{j+1}^{(\theta)} \ebold[k,\theta^{\diamond j}\star k]\right)\otimes a^{(\theta)}.
\end{equation}
constructed by using these data. 
Assume that
\begin{eqnarray}\label{equation:InvertibilityHyp}
&&\mbox{$G=A^{-1}\in \Mfrak((1/z))$ exists, and}\\
\label{equation:gUpperLeft}
&&g(i,j)=G(i-1,j-1)\;\;\mbox{for $i,j=1,\dots,n$.}\end{eqnarray}
Then  \eqref{equation:GSDE2}, \eqref{equation:GSDE1}, and \eqref{equation:GSDE3} hold, i.e., the data \eqref{equation:Data} constitute a solution of the generalized Schwinger-Dyson equation.
\end{Proposition}
 We complete the proof below in \S\ref{subsection:ClosingArgument}
after deducing Proposition \ref{Proposition:MainResultBlockManipulation}
from  Proposition \ref{Proposition:GSDE}
and proving a  lemma excusing us from having to verify \eqref{equation:GSDE3}
when proving Proposition \ref{Proposition:GSDE}.
The minus signs inserted in definition \eqref{equation:Adef}
turn out to save us from writing many minus signs later.
\subsubsection{Remark}
We can rewrite \eqref{equation:Adef} as
\begin{equation}\label{equation:AdefBis}
A=-\one\otimes a^{(0)}-\sum_{\theta=1}^q\left(
\hat{\lambda}^{(\theta)}+\sum_{j=1}^\infty\kappa_{j+1}^{(\theta)} \lambda^{(\theta)}\right)\otimes a^{(\theta)}
\end{equation}
in terms of the lowering and raising operators considered in Remark \ref{Remark:LoweringAndRaising},
thus making the relationship of the definition of $A$ to the setup of \cite{Voiculescu} more transparent.

\subsubsection{Completion of the proof of Proposition \ref{Proposition:MainResultBlockManipulation}
with Proposition \ref{Proposition:GSDE} granted}
We identify $X$ in Proposition \ref{Proposition:MainResultBlockManipulation}
with $X$ in Proposition \ref{Proposition:DigitalLinearization}.
We complete the choice of positive integer $q$ and the given family
$\left\{\left\{\kappa_j^{(\theta)}\right\}_{\theta=1}^q\right\}_{j=2}^\infty$
to a family
\begin{equation}\label{equation:LinearData}
\left(q,n,g,a^{(0)},\left\{a^{(\theta)}\right\}_{\theta=1}^q,\left\{\left\{\kappa_j^{(\theta)}\right\}_{j=2}^\infty\right\}_{\theta=1}^q\right)
\end{equation}
of the form \eqref{equation:Data} where
\begin{eqnarray*}
n&=&p+N>p,\\
a^{(0)}&=&-\left(L_0+\left[\begin{array}{cc}
zI_p&0\\
0&0\end{array}\right]-\sum_{\theta=1}^qL_\theta\right)\in \Mat_n(\CC[z]),\\
a^{(\theta)}&=&-L^{(\theta)}\in \Mat_n(\CC)\;\;\mbox{for $\theta=1,\dots,q$ and}\\
g(i,j)&=&-L^{-1}(i-1,j-1)\;\;\mbox{for $i,j=1,\dots,n$.}
\end{eqnarray*}
For the family \eqref{equation:LinearData}
the hypotheses \eqref{equation:InvertibilityHyp} and \eqref{equation:gUpperLeft}
of Proposition \ref{Proposition:GSDE} are fulfilled by \eqref{equation:Ldef} and
the definition of $g$, respectively.
Thus \eqref{equation:LinearData} is a solution of the generalized Schwinger-Dyson equation,
i.e., \eqref{equation:GSDE2}, \eqref{equation:GSDE1},  and \eqref{equation:GSDE3} hold.
Property \eqref{equation:GSDE5} holds by construction and
property \eqref{equation:GSDE0} holds by \eqref{equation:GeneralizedResolventKey}.
The proof of Proposition \ref{Proposition:MainResultBlockManipulation} is complete
modulo the proof of Proposition \ref{Proposition:GSDE}.
\qed

\begin{Lemma}\label{Lemma:Bootstrap} To prove Proposition \ref{Proposition:GSDE} it is necessary only to
verify statements \eqref{equation:GSDE2} and \eqref{equation:GSDE1}
for  data \eqref{equation:Data} satisfying hypotheses \eqref{equation:InvertibilityHyp}
and \eqref{equation:gUpperLeft}.
\end{Lemma}
 In the proof below we are reusing elements of the  ``secondary trick'' used in \cite{Anderson} to obtain certain correction terms. 
 \proof
The weakened version of
Proposition \ref{Proposition:GSDE} delivering only conclusions  \eqref{equation:GSDE2} 
and \eqref{equation:GSDE1} for data \eqref{equation:Data} satisfying \eqref{equation:InvertibilityHyp} and \eqref{equation:gUpperLeft} we will call Proposition \ref{Proposition:GSDE}$-\epsilon$. Our task is to derive Proposition \ref{Proposition:GSDE} from Proposition \ref{Proposition:GSDE}$-\epsilon$. To that end fix \linebreak
$b\in \Mat_n(\CC((1/z))$ arbitrarily
and consider new data consisting of
\begin{equation}\label{equation:NewData} 
\left\{\begin{array}{l}
\mbox{a positive integer $\hat{n}=2n$ (but $q$ the same as before),}\\
\mbox{a matrix $\hat{a}^{(0)}=\left[\begin{array}{cc}
a^{(0)}&b\\
0&a^{(0)}
\end{array}\right]\in \Mat_{\hat{n}}(\CC((1/z)))$,}\\
\mbox{matrices $\hat{a}^{(\theta)}=\left[\begin{array}{cc}
a^{(\theta)}&0\\
0&a^{(\theta)}\end{array}\right]\in \Mat_{\hat{n}}(\CC)$ for $\theta=1,\dots,q$,}\\
\mbox{a matrix $\hat{g}=\left[\begin{array}{cc}
g&h\\
0&g
\end{array}\right]\in \Mat_{\hat{n}}(\CC((1/z)))$ ($h$ to be determined), and}\\
\mbox{a  family $\{\{\kappa_j^{(\theta)}\}_{j=2}^\infty\}_{\theta=1}^q$  of complex numbers (same as before).}
\end{array}\right.
\end{equation}
We will apply Proposition \ref{Proposition:GSDE}$-\epsilon$ to the new data \eqref{equation:NewData}
thereby deriving \eqref{equation:GSDE3} for the old data \eqref{equation:Data}.
To apply Proposition \ref{Proposition:GSDE}$-\epsilon$ we need first to verify invertibility of the matrix
$$
\hat{A}=-\hat{a}^{(0)}-\sum_{\theta=1}^q\sum_{k\in \NN}
\left(\ebold[\theta\star k,k]+\sum_{j=1}^\infty\kappa_{j+1}^{(\theta)} \ebold[k,\theta^{\diamond j}\star k]
\right)\otimes \hat{a}^{(\theta)}.
$$
Because $\natural$ is an isomorphism,
there exists unique $\tilde{A}\in \Mat_n(\Mfrak((1/z)))$ such that $(\tilde{A})^\natural=A$.
Using \eqref{equation:Switcharoo} we obtain the relation
$$\hat{A}=\left[\begin{array}{cc}
\tilde{A}&- b\otimes \one\\
0&\tilde{A}
\end{array}\right]^\natural\in \Mfrak((1/z)),$$
and we have explicitly
$$\hat{G}=\left(\left[\begin{array}{cc}
\tilde{A}&-b\otimes \one\\
0&\tilde{A}
\end{array}\right]^{-1}\right)^\natural=\left[\begin{array}{cc}
\tilde{A}^{-1}&\tilde{A}^{-1}(b\otimes \one)\tilde{A}^{-1}\\
0&\tilde{A}^{-1}
\end{array}\right]^\natural.$$
Thus  the new data \eqref{equation:NewData} satisfy \eqref{equation:InvertibilityHyp},
and moreover there is a choice of $h$ we can in principle read off from the last displayed line above
so that hypothesis \eqref{equation:gUpperLeft} is satisfied. (An explicit formula for $h$ is not needed.)
Statement \eqref{equation:GSDE1} of Proposition \ref{Proposition:GSDE}$-\epsilon$
applied to the new data asserts that
$$\lim_{j\rightarrow\infty}\val\,\left[\begin{array}{cc}
a^{(\theta)} g&a^{(\theta)} h\\
0&a^{(\theta)} g
\end{array}\right]^j=-\infty\;\;\mbox{for $\theta=1,\dots,q$.}$$
This can also be deduced directly from \eqref{equation:GSDE1} as it pertains to the old data \eqref{equation:Data}.
Finally, the key point is that by statement \eqref{equation:GSDE2} as it pertains to the new data
\eqref{equation:NewData} we have
$$I_{\hat{n}}+\left[\begin{array}{cc}
a^{(0)}&b\\
0&a^{(0)}
\end{array}\right]\left[\begin{array}{cc}
g&h\\
0&g\end{array}\right]+\sum_{\theta=1}^q\sum_{j=2}^\infty \kappa_j^{(\theta)}\left(\left[\begin{array}{cc}
a^{(\theta)}&0\\
0&a^{(\theta)}
\end{array}\right]\left[\begin{array}{cc}
g&h\\
0&g
\end{array}\right]\right)^j=0.
$$
Looking in the upper left corners, we obtain an identity
$$a^{(0)}h+b g+\sum_{\theta=1}^q\sum_{j=2}^\infty
\sum_{\nu=0}^{j-1}\kappa_j^{(\theta)}(a^{(\theta)}g)^\nu(a^{(0)}h)(a^{(\theta)}g)^{j-1-\nu}=0.
$$
The latter equation, because $b$ is arbitrary and $g$ is invertible
by \eqref{equation:GSDE2} as it pertains to the old data \eqref{equation:Data}, proves  that
\eqref{equation:GSDE3} holds for the old data. In other words,  Proposition \ref{Proposition:GSDE}$-\epsilon$ does indeed imply
Proposition \ref{Proposition:GSDE}.  \qed

\subsection{Proof of Proposition \ref{Proposition:GSDE}}\label{subsection:ClosingArgument}
In broad outline the proof is
similar to the proof we previously gave for Lemma \ref{Lemma:HessenbergToeplitz}.
But more machinery is needed.
\subsubsection{Block decompositions}
Throughout the proof it will be convenient to work
with the block decompositions defined by the formulas
 $$A=\sum_{i,j\in \NN}\ebold[i,j]\otimes A\langle i,j\rangle\;\;
 \mbox{and}\;\;G=\sum_{i,j\in \NN}\ebold[i,j]\otimes G\langle i,j\rangle$$
 where
$$A\langle i,j\rangle, G\langle i,j\rangle\in \Mat_n(\CC((1/z)))\;\;
\mbox{and in particular}\;\;G\langle 0,0\rangle=g.$$

The key point is contained in the following result which says that the first column of blocks in $G$
has a relatively simple structure.
\begin{Lemma}\label{Lemma:Exploit}
We have
\begin{equation}\label{equation:LastRecursion}
G\langle\theta_1\star \cdots\star \theta_k,0\rangle=ga^{(\theta_1)}ga^{(\theta_2)}\cdots  ga^{(\theta_k)}g
\end{equation}
for $k\in\NN$ and $\theta_1,\dots,\theta_k\in \{1,\dots,q\}$.
\end{Lemma}
\noindent This relation generalizes formula \eqref{equation:HessTopClincherBis} above.
\proof For the proof we will use Lemma \ref{Lemma:AbstractSchurComplement} in the case
\begin{eqnarray*}
&&\AAA=\Mfrak((1/z)),\;\;A\in \AAA:\;\mbox{as on line \eqref{equation:Adef}},\\
&&\pi=\ebold[0,0]\otimes I_n\in \AAA\;\;\mbox{and}\;\;\sigma=\sum_{i=1}^\infty \ebold[i,i]\otimes I_n\in \AAA.
\end{eqnarray*}
We will also use the matrices
\begin{eqnarray*}
&&R^{(\theta)}=\sum_{k\in \NN}\ebold[k\star \theta,k ]\otimes I_n\in \Mfrak\;\;\mbox{and}\;\;
\hat{R}^{(\theta)}=\sum_{k\in \NN}\ebold[k,k\star \theta]\otimes I_n\in \Mfrak\\
\nonumber&&\hspace{3.0in}\mbox{for $\theta=1,\dots,q$.}
\end{eqnarray*}
These matrices satisfy
\begin{eqnarray}
\label{equation:RIdentity}
&&\hat{R}^{(\theta)}R^{(\theta')}=\delta_{\theta\theta'}\one\;\;\mbox{for $\theta,\theta'=1,\dots,q$,}\\
\label{equation:RARidentity}
&&\hat{R}^{(\theta)}AR^{(\theta')}=\delta_{\theta\theta'}A\;\;\mbox{for $\theta,\theta'=1,\dots,q$, and}\\
\label{equation:RIdentityBis}
&&\sum_{\theta=1}^q R^{(\theta)}\hat{R}^{(\theta)}=\sigma,
\end{eqnarray}
as can either be verified directly by straightforward if tedious calculation or by using the 
definitions and relations \eqref{equation:Sinister}, \eqref{equation:ExcellentFamilyBis}, \eqref{equation:Dextrous},
\eqref{equation:Commutation1}, and \eqref{equation:Commutation2} listed above in Remark \ref{Remark:LoweringAndRaising}, along with the rewrite \eqref{equation:AdefBis}
of the definition of $A$.
It follows by \eqref{equation:RARidentity} and  \eqref{equation:RIdentityBis} that
$$
\sigma A\sigma=\left(\sum_{\theta=1}^q R^{(\theta)}\hat{R}^{(\theta)}\right)A\left(\sum_{\theta=1}^q R^{(\theta')}\hat{R}^{(\theta')}\right)
=\sum_{\theta=1}^q R^{(\theta)}A\hat{R}^{(\theta)}.
$$
It follows in turn that $A_\sigma^{-1}$ exists and more precisely that
\begin{equation}\label{equation:LeanIn}
A_\sigma^{-1}=\sum_{\theta=1}^q R^{(\theta)}G\hat{R}^{(\theta)},
\end{equation}
as one verifies using \eqref{equation:RIdentity}.
We furthermore have
$$\sigma A\pi A^{-1}\pi=-\sum_{\theta=1}^q (\ebold[\theta,0]\otimes a^{(\theta)})(\ebold[0,0]\otimes g)=
-\sum_{\theta=1}^q R^{(\theta)}(\ebold[0,0]\otimes a^{(\theta)}g),$$
as one can immediately check.
Finally, we have the following chain of equalities:
\begin{eqnarray*}
\sum_{\theta=1}^q\sum_{i\in \NN} \ebold[i\star \theta,0]\otimes G\langle i\star\theta,0\rangle&=&\sigma A^{-1}\pi\;=\;-A_\sigma^{-1}A\pi A^{-1}\pi\\
&=&\sum_{\theta=1}^q \sum_{\theta'=1}^qR^{(\theta)}G\hat{R}^{(\theta)}R^{(\theta')}(\ebold[0,0]\otimes a^{(\theta')}g)\\
&=&\sum_{\theta=1}^q\sum_{i\in \NN} \ebold[i\star \theta,0]\otimes G\langle i,0\rangle a^{(\theta)} g.
\end{eqnarray*}
At the first step we used Lemma \ref{Lemma:Unpin}, 
at the second equation \eqref{equation:SC3} of Lemma \ref{Lemma:AbstractSchurComplement},
at the third \eqref{equation:LeanIn}, and
at the last  \eqref{equation:RIdentity}.
Thus \eqref{equation:LastRecursion} holds.
\qed

\subsubsection{Proof of  \eqref{equation:GSDE1}}
By definition of $\Mfrak((1/z))$ we have
\begin{equation*}
\lim_{i\rightarrow\infty}\val\, G\langle i,0\rangle=-\infty.
\end{equation*}
Thus by Lemma \ref{Lemma:Exploit} we have for $\theta=1,\dots,q$ that
$$\lim_{i\rightarrow\infty}\val\,(a^{(\theta)}g)^i=\lim_{i\rightarrow\infty}\val\, a^{(\theta)}G\langle \theta^{\diamond i},0\rangle=0,$$
which proves statement \eqref{equation:GSDE1}.

\subsubsection{Proof of \eqref{equation:GSDE2}} \label{subsubsection:ClosingArgument}
Consider the following calculation:
\begin{eqnarray*}
I_n&=&\sum_{j\in \NN}A\langle 0,j\rangle G\langle j,0\rangle\;=\;- a^{(0)} G\langle 0,0\rangle-\sum_{\theta=1}^q\sum_{j=1}^\infty \kappa_{j+1}^{(\theta)}a^{(\theta)} G\langle \theta^{\diamond j},0\rangle\\
&=&- a^{(0)} g-\sum_{\theta=1}^q\sum_{j=1}^\infty \kappa_{j+1}^{(\theta)}(a^{(\theta)}g)^j.
\end{eqnarray*}
The first step holds by definition of $G$, the second by definition of $A$,
and the last by Lemma \ref{Lemma:Exploit}.
This calculation proves statement \eqref{equation:GSDE2}.

\subsubsection{Completion of the proof}
By Lemma \ref{Lemma:Bootstrap} it is necessarily the case that statement \eqref{equation:GSDE3} holds.
The proof of Proposition \ref{Proposition:GSDE} is complete and in turn the proof of Proposition \ref{Proposition:MainResultBlockManipulation} is complete. \qed
\subsection{Miscellaneous remarks}
\subsubsection{}
In Lemma \ref{Lemma:Exploit} and its proof we are exploiting
the type of recursion used, for example, in \cite{NagWoess}, and used more generally in many investigations of random walk on infinite trees.

\subsubsection{}
By exploiting the identities \eqref{equation:RIdentity}, \eqref{equation:RARidentity},
and \eqref{equation:RIdentityBis}, we are reusing some features of the proof of \cite[Prop. 9]{Anderson}.

\subsubsection{}
In the recently introduced bi-free probability setup \cite{VoiculescuJanus}
the left and right variants of lowering and raising operators play an equal role.
Our approach here is potentially useful for calculating the laws of matrix-polynomial combinations
of bi-free collections of noncommutative random variables.

\subsubsection{}
By continuing the line of argument in the proof of Lemma \ref{Lemma:Exploit}
and using \eqref{equation:SC4} of Lemma \ref{Lemma:AbstractSchurComplement},
it is not difficult to prove  that
\begin{equation}\label{equation:WholeEnchilada}
GR^{(\theta)}=(R^{(\theta)}-G\pi AR^{(\theta)})G\;\;\mbox{for $\theta=1,\dots,q$.}
\end{equation}
Since statement \eqref{equation:WholeEnchilada} is not needed for the proof of Theorem \ref{Theorem:MainResult},
we omit its proof. It is easy to see that Lemma \ref{Lemma:Exploit} and  \eqref{equation:WholeEnchilada}  together allow one to make
every block $G\langle i,j\rangle$ explicit in terms of $g$ and $A$.  Doing so in a systematic way
one would obtain a generalization of Theorem \ref{Theorem:MainResult} 
having the full statement of Theorem \ref{Theorem:Aomoto} as a consequence.

\section{Notes on Newton-Puiseux series}\label{section:NewtonPolygon}
At this point in the paper we switch from the noncommutative viewpoint
previously stressed to the viewpoint of commutative algebra and  algebraic geometry.
The latter is maintained throughout the rest of the paper.

 \subsection{Newton-Puiseux series and Newton polygons}\label{subsection:Newton}
We review  basic devices for understanding singularities of plane algebraic curves in characteristic zero
and make some definitions needed for later calculations.

\subsubsection{The algebraic closure of $\CC((1/z))$}
Let $\KK$ denote the union of the tower of fields $\{\CC((1/z^{1/n!}))\}_{n=1}^\infty$. In other words, $\KK$ arises by adjoining roots of $z$ 
of all orders to $\CC((1/z))$.
We call an element of $\KK$ a {\em Newton-Puiseux series}.
When discussing $\KK$ below we often use the more apposite abbreviated notation 
$\KK_0=\CC((1/z))$.
It is well-known that $\KK$ is the algebraic closure of $\KK_0$. 
See, e.g.,  \cite[Cor. 13.15]{Eisenbud}. The original insight is due to  Newton. 

\subsubsection{Extension of the valuation function $\val$ to $\KK$}
Each element $f\in \KK$ has  by definition a unique {\em Newton-Puiseux expansion}
$f=\sum_{u\in \QQ}c_uz^u$
with coefficients $c_u\in \CC$ such that for some positive integer $N$ depending on $f$ one has
$c_u=0$ unless $u\leq N$ and $Nu\in \ZZ$. To extend to $\KK$ the valuation defined on $\KK_0$,
we define
$$\val\, f=\sup \{u\in \QQ\mid c_u\neq 0\}\in \QQ\cup\{-\infty\}.$$
The properties \eqref{equation:VP1}, \eqref{equation:VP2}, and \eqref{equation:VP3}
continue to hold for the extension of $\val$ to $\KK$. 
\begin{Proposition}\label{Proposition:NewtonPolygon}
Let $P(y)\in \KK_0[y]$ be a polynomial of degree $n$ in a variable $y$ with coefficients in the field $\KK_0$. 
Write 
$$P(y)=\sum_{i=0}^{n}a_i y^{n-i}=a_0\prod_{i=1}^n (y-r_i)\;\;
(a_i\in \KK_0\;\mbox{and}\;r_i\in \KK),$$ 
enumerating the roots $r_i$ so that $\val\,r_1\geq \cdots \geq \val\,r_n$.
(i) Then we have
\begin{eqnarray}\label{equation:NewtonPolygon}
&&\val\, a_i\leq \val\, a_0+\val\,r_1\cdots r_i\;\mbox{for $i=0,\dots,n$, }\\
\label{equation:NewtonPolygonBis}
&&\mbox{with equality for  $i=0,n$
and furthermore}\\
\nonumber&&\mbox{for $i=1,\dots,n-1$ s.t. $\val\,r_i>\val\,r_{i+1}$.}
\end{eqnarray}
(ii)  Let $\psi:[0,n]\rightarrow\RR\cup\{-\infty\}$ 
be the infimum of all affine linear functions \linebreak $\lambda:[0,n]\rightarrow\RR$
satisfying $\lambda(i)\geq \val\, a_i$ for $i=0,\dots,n$.  Then 
we have the integral formula $\psi(s)=\val\, a_0+
\int_0^s \val\,r_{\lceil u\rceil}\,du$.
\end{Proposition}
Here $\lceil u\rceil$ denotes the least integer not less than $u$. The function $\psi$ is concave by construction.
The function $\psi$  (or rather, its graph)
is the {\em Newton polygon} 
associated with $P(y)$,
up to reflections in 
and translations parallel to the horizontal and vertical axes
(conventions vary from author to author).
 Similarly, a Newton polygon is
attached to any one-variable polynomial
with coefficients in a discretely valued field.
For background on Newton polygons see \cite[Chap. 2, Sec. 5]{Artin} or 
\linebreak \cite[Part III, Chap. 8, Sec. 3]{Brieskorn}. 

\proof (i) Since $(-1)^ia_{i}/a_0$ for $i>0$ is the $i^{th}$ symmetric function of the roots $r_1,\dots,r_n$, the result follows straightforwardly
from \eqref{equation:VP1}, \eqref{equation:VP2} and \eqref{equation:VP3}. (ii) 
The deduction of this statement from the preceding one is standard. 
\qed

\subsubsection{Extension of the valuation to $\Mat_n(\KK)$}
We extend $\val\,$ from $\Mat_n(\KK_0)$
to  $\Mat_n(\KK)$ by the rule $\val\,A=\max_{i,j=1}^n \val\,A(i,j)$. Then $\Mat_n(\KK)$ satisfies all the axioms
of a Banach algebra over $\KK$ except completeness.
Lack of completeness will not be an issue.

\subsection{Applications}

We present several applications of the preceding machinery needed
for the proof of Proposition \ref{Proposition:MainResultAlgebraicityProof}.

\subsubsection{Specialized matrix notation}\label{subsubsection:eNotation}
Given $A\in \Mat_n(\CC((1/z)))$, we define
$$e(A)=(e_1(A),\dots,e_n(A))\in \CC((1/z))^n$$
by the formula
$$\det(tI_n-A)=t^n+\sum_{i=1}^n e_i(A)(-1)^i t^{n-i}\in \CC((1/z))[t].$$
The Cayley-Hamilton Theorem takes then the form
\begin{equation}\label{equation:CayleyHamilton}
A^n+\sum_{i=1}^n (-1)^ie_i(A)A^{n-i}=0.
\end{equation}

\begin{Proposition}\label{Proposition:NegativeSpectralValuation}
For  $A\in \Mat_n(\CC((1/z)))$ the following are equivalent:
\begin{enumerate}
\item[(I)]  $e(A)\in (1/z)\CC[[1/z]]^n$.
\item[(II)] $\lim_{k\rightarrow\infty} \val\,A^{k}=-\infty$.
\end{enumerate}
\end{Proposition}
\proof We add one further statement to the list above:
\begin{enumerate}
\item[(III)] Every eigenvalue of $A$ in $\KK$ has (strictly) negative valuation.
\end{enumerate}
Statements (I) and (III) are equivalent by Proposition \ref{Proposition:NewtonPolygon}.
It remains only to prove the equivalence (II)$\Leftrightarrow$(III). It is actually easier to prove more.
We will prove the equivalence (II)$\Leftrightarrow$(III) for $A\in \Mat_n(\KK)$.
Supposing at first that $A$ consists of a single Jordan block, one verifies the equivalence
by inspection. In general we can write $A=WJW^{-1}$
where $W\in \GL_n(\KK)$ and $J\in \Mat_n(\KK)$ is block-diagonal with diagonal blocks of the Jordan form
and we have a bound
$$\left|\val\, A^k -\val\, J^k\right|\leq \val\, W^{-1} +\val\, W$$
which establishes the equivalence (II)$\Leftrightarrow$(III) in general.
\qed

\subsubsection{Negative spectral valuation}
We say that $A\in \Mat_n(\KK_0)$ has {\em negative spectral valuation}
if the equivalent conditions (I) and (II) above hold.

\subsubsection{Algebraic and nonsingular algebraic elements of $\CC[[t]]$}
\label{subsubsection:NonsingularAlgebraic}
In some situations $1/z$ rather than $z$ is the natural parameter to work with.
We therefore make the following definitions which are in principle redundant but in practice convenient.
We say that $f(t)\in \CC[[t]]$ is {\em algebraic}
 if  $F(t,f(t))=0$ for some $0\neq F(x,y)\in \CC[x,y]$.
Of course $f(t)\in \CC[[t]]$ is algebraic if and only if $f(1/z)\in  \CC((1/z))$ is algebraic in the sense  defined in \S\ref{subsubsection:Algebraic}.
 We call $f(t)$ {\em nonsingular} algebraic if
there exists $F(x,y)\in \CC[x,y]$  satisfying $\frac{\partial F}{\partial y}(0,f(0))\neq 0$ and $F(t,f(t))=0$.

The next statement is the key to desingularization.
\begin{Lemma}\label{Lemma:NewtonPolygonBis}
Let $f(t)\in t\CC[[t]]$ be a power series with vanishing constant term.
Let $F(x,y)\in\CC[x,y]$ be a polynomial not divisible by $x$ such that $F(t,f(t))=0$.
Then the following conditions are equivalent:
\begin{enumerate}
\item[(I)] $F(1/z,y)\in \KK_0[y]$ has exactly one root in $\KK$ of negative valuation.
\item[(II)] $\frac{\partial F}{\partial y}(0,0)\neq 0$.
\end{enumerate}
\end{Lemma}
\proof Write $F(x,y)=\sum_{i=0}^n p_i(x)y^{n-i}$ where $n>0$, $p_i(x)\in \CC[x]$ and $p_0(x)\neq 0$.
For $p(x)\in \CC[x]$, let $\ord\,p(x)$ denote the exponent of the highest power of $x$ dividing $p(x)$.
Let $\phi:[0,n]\rightarrow\RR\cup\{+\infty\}$ denote the supremum of all affine linear functions
$\lambda:[0,n]\rightarrow\RR$ such that $\lambda(i)\leq \ord\,p_i(x)$ for $i=0,\dots,n$.
The function $\phi$ is convex. Since $F(0,0)=0$ we have $\phi(n)>0$.
Since $x$ does not divide $F(x,y)$, we have $\min_{i=0}^n \ord\, p_i(x)=\min_{i=0}^n \phi(i)=0$.
Let $i_0$ be the maximum of $i=0,\dots,n-1$ such that $\phi(i_0)=0$.
Statement (ii) is equivalent to the assertion that $i_0=n-1$.
Now $F(1/z,y)$ has $n-i_0$ roots of negative valuation by Proposition \ref{Proposition:NewtonPolygon}
and the observation that $\ord\, p(x)=-\val\, p(1/z)$.
Thus statement (I) is equivalent to the assertion that $i_0=n-1$.
\qed

The next statement summarizes just enough of the theory of resolution of singularities of plane algebraic curves
in characteristic zero for our purposes. 

\begin{Proposition} \label{Proposition:Desingularization}
Let $\sum_{i=0}^\infty c_it^i\in t\CC[[t]]\;\;(c_i\in \CC)$
  be an algebraic power series. Then 
$\sum_{i=N}^\infty c_{i+N}t^i\in t^N\CC[[t]]$ is a nonsingular algebraic power series for all $N\gg 0$.
\end{Proposition}
\proof 
Let $f=\sum_{i=1}^\infty c_iz^{-i}\in \CC((1/z))$ and $f_N=\sum_{i=N}^\infty c_{i+N}z^{-i}\in \CC((1/z))$. 
Let $F(x,y)\in \CC[x,y]$ (resp., $F_N(x,y)\in \CC[x,y]$)
denote the irreducible equation (see \S\ref{subsubsection:Algebraic}) of $f$ (resp., $f_N$).  It will be enough to show that $F_N(1/z,y)\in \KK_0[y]$ has exactly one root in $\KK$
of negative valuation for $N\gg 0$.
If $f_{N_0}=0$ for some $N_0$ then $f_N=0$
for all $N\geq N_0$
and there is nothing to prove.
Thus we may assume without loss of generality that $f_N\neq 0$ for all $N\geq 0$.
Let $n$ denote  the dimension of $\CC(z,f)$ over $\CC(z)$.
It is clear that $f$ and $f_N$ generate the same extension of $\CC(z)$. 
Thus $F(x,y)$ and $F_N(x,y)$ have the same degree in $y$, namely $n$.
 Let $r_1,\dots,r_n$ denote the roots in $\KK$ of the polynomial $F(1/z,y)\in \KK_0[y]$,
enumerated so that $r_n=f$.  Let  $h_N=\sum_{i=0}^{2N-1}a_iz^{-i}=f-f_N/z^{N}\in \CC(z)$, in which case necessarily $f_N=z^N(f-h_N)$.
Then for a suitable enumeration $r_{1,N},\dots,r_{n,N}$ of the roots in $\KK$ of the polynomial $F_N(1/z,y)\in \KK_0[y]$,
we have \linebreak
$r_{i,N}=z^N(r_i-h_N)$ for $i=1,\dots,n$ and $r_{n,N}=f_N$.
Now the roots $r_1,\dots,r_n$ are distinct due to irreducibility of $F(x,y)$,
and clearly $\val\,(f-h_N)\leq-2N$. 
Because $h_N\rightarrow_{N\rightarrow\infty}f$ with respect to the valuation $\val$,
it follows that  for some integer $N_0>0$ depending only on $f$, and all integers $N\geq N_0$, we have
$$
\min_{i=1}^{n-1}\val\,r_{i,N}=N+\min_{i=1}^{n-1}\val(r_i-f)\geq 0>-N\geq \val\, f_N=\val\,r_{n,N},
$$
whence the result via Lemma \ref{Lemma:NewtonPolygonBis}. \qed

\section{Evaluation of algebraic power series on matrices}\label{section:Gizmo}
The main result of this section is Proposition \ref{Proposition:Gizmo} below.
The main tools used in this section are the Cayley-Hamilton Theorem, the Weierstrass Preparation Theorem, 
and Proposition \ref{Proposition:Desingularization} above.

\subsection{Motivation}
Let $A$ be an $n$-by-$n$ matrix with complex entries.
We take our inspiration
from the undergraduate level approach in \cite{WilliamsonTrotter} to computing the matrix exponential $\exp(tA)$.
The approach is lengthy as presented for sophomores but it can be summarized quickly at graduate level as follows.
Perform Weierstrass division (possible globally in this case) in order to obtain an identity relating two-variable entire functions of complex variables $t$ and $X$, namely
\begin{equation}\label{equation:UrWeierstrass}
\exp(tX)=\sum_{k=0}^{n-1}y_k(t)X^k+Q(X,t)\det(XI_n-A),
\end{equation}
for  suitable and unique remainder $\sum_{k=0}^{n-1}y_k(t)X^k$ and quotient $Q(X,t)$.
Differentiation of \eqref{equation:UrWeierstrass} on both sides with respect to $t$ 
yields a first order homogeneous linear differential equation for the vector function
$y(t)=(y_1(t),\dots,y_k(t))$
which together with the evident initial value data uniquely determines $y(t)$.
One can then go on to solve explicitly for the functions $y_k(t)$ in closed form.
By plugging in $X=A$ on both sides of \eqref{equation:UrWeierstrass}
and using the Cayley-Hamilton Theorem one then has finally
\begin{equation}\label{equation:UrWeierstrassBis}
\exp(tA)= \sum_{k=0}^{n-1}y_k(t)A^k.
\end{equation}
Formula \eqref{equation:UrWeierstrassBis}
makes no reference to the Jordan canonical form of $A$.
Indeed, by construction, the coefficients $y_k(t)$ are uniquely
determined by the characteristic polynomial of $A$ alone.
We will make roughly analogous use of Weierstrass division below to evaluate
algebraic power series on matrices with entries in $\CC((1/z))$ of negative spectral valuation. 

\subsection{$I$-adic convergence, power series,  and Weierstrass division}\label{subsection:adic}
We pause to review generalities connected with the formal power series
version of the Weierstrass Preparation Theorem.

\subsubsection{$I$-adic convergence}
Given a commutative ring $R$ with unit, an  ideal $I$,
and a sequence $\{a\}\cup\{a_i\}_{i=1}^\infty$ in $R$, one says $\lim_{i\rightarrow\infty}a_i=a$
holds {\em $I$-adically} if for every positive integer $k$ there exists a positive integer $i_0=i_0(k)$
such that \linebreak $a-a_i\in I^k$ for all $i\geq i_0$.  Similarly, one can speak of $I$-adic Cauchy sequences and $I$-adic completeness. Consider, e.g., the ring $\CC[[u_1,\dots,u_n]]=\CC[[u]]$
and the maximal ideal $I=(u_1,\dots,u_n)\subset \CC[[u]]$.
Then $f_i\in \CC[[u]]$ converges $I$-adically to $f\in \CC[[u]]$
if and only if for every $n$-tuple $(\nu_1,\dots,\nu_n)$ of nonnegative integers
and every sufficiently large index $i$ depending on $(\nu_1,\dots,\nu_n)$,
the Taylor coefficient
 $\frac{1}{\nu_1!\cdots \nu_n!}\frac{\partial^{\nu_1+\cdots+\nu_n}f_i}{\partial u_1^{\nu_1}\cdots \partial u_n^{\nu_n}}(0)$ equals the Taylor coefficient $\frac{1}{\nu_1!\cdots \nu_n!}\frac{\partial^{\nu_1+\cdots+\nu_n}f}{\partial u_1^{\nu_1}\cdots \partial u_n^{\nu_n}}(0)$. It is easy to see that the ring $\CC[[u]]$ is $I$-adically complete.

\subsubsection{Weierstrass division}\label{subsubsection:Weierstrass}
We now briefly recall the {\em Weierstrass Preparation Theorem} from a more active point of view
emphasizing the algorithm of Weierstrass division.
See, e.g., \cite[Thm. 5, p. 139, Chap. VII, \S 1]{ZariskiSamuel} for background and proof.
The theorem concerns an $(n+1)$-variable power series ring over a field, with one of the variables
singled out for special treatment.
For definiteness  we take the coefficient field to be $\CC$. Consider the ring $\CC[[u_1,\dots,u_n,t]]=\CC[[u,t]]$,
with $t$ distinguished.
One is given a {\em divisand} $F(u,t)\in \CC[[u,t]]$
and a {\em divisor} $D(u,t)\in \CC[[u,t]]$. Of the latter it is assumed that there exists a positive integer $m$
(called the {\em multiplicity} of the divisor) such that $D(0,t)=t^mU(t)$ for some $U(t)\in \CC[[t]]$ such that $U(0)\neq 0$.
The Weierstrass division process delivers a  {\em quotient} $Q(u,t)\in \CC[[u,t]]$
and a {\em remainder} $R(u,t)\in \CC[[u]][t]$. The pair $(Q(u,t),R(u,t))$
is uniquely determined by two requirements.
Firstly,
the division equation  $F(u,t)=Q(u,t)D(u,t)+R(u,t)$ must hold.
Secondly, $R(u,t)$ must be a polynomial in $t$ of degree $<m$.  It bears emphasis that
if $D(u,t)\in \CC[[u]][t]$ is monic of degree $m$ such that $D(0,t)=t^m$,
and $F(u,t)\in \CC[[u]][t]$, then the Euclidean (i.e., high school) and Weierstrass division processes 
deliver the same quotient and remainder.

\begin{Lemma}\label{Lemma:Approximation}
We continue in the setting of the preceding paragraph.
However, for simplicity we assume now 
that $D(u,t)$ is monic of degree $m$ such that \linebreak $D(0,t)=t^m$.
Consider the ideal $I=(u_1,\dots,u_n)\subset \CC[[u]][t]$.
Let $k$ be a positive integer.
If $t^k$ divides $F(u,t)$, then $R(u,t)$ belongs to
the ideal $I^{\lfloor k/m\rfloor}$.
(Here $\lfloor c\rfloor$ denotes the greatest integer not exceeding $c$.)
\end{Lemma}
 It follows that formation of Weierstrass remainder upon division by $D(u,t)$ viewed as a 
function from $\CC[[u,t]]$ to $\CC[[u]][t]$ is continuous with respect to the $(t)$-adic topology
on the source and the $(u_1,\dots,u_n)$-adic topology on the target.
\proof Let $F_0(u,t)=F(u,t)/t^k$. Let $R_0(u,t)$ denote the remainder
of $F_0(u,t)$ upon Weierstrass division by $D(u,t)$.
Then $R(u,t)$ is the remainder of $t^kR_0(u,t)$ upon high school
division by $D(u,t)$. This noted, there is no loss of generality
in assuming that $F(u,t)=t^k$. 
Write $D(u,t)=t^m+\sum_{i=0}^{m-1} a_it^i$ with coefficients \linebreak $a_i=a_i(u)\in \CC[[u]]$
such that $a_i(0)=0$. Write  $R(u,t)=\sum_{i=0}^{m-1}b_it^i$ with coefficients
$b_i=b_i(u)\in \CC[[u]]$.
Then we have
$$
\left[\begin{array}{ccccc}
&&&-a_0\\
1&&&\vdots\\
&\ddots&&\vdots\\
&&1&-a_{m-1}
\end{array}\right]^k
\left[\begin{array}{c}
1\\
0\\
\vdots\\
0\end{array}\right]=
\left[\begin{array}{c}
b_0\\
\vdots\\
b_{m-1}
\end{array}\right],$$
where the matrix on the left is the so-called {\em companion matrix}
for $D(u,t)$.
Clearly every entry of the $m^{th}$ power of the companion matrix belongs
to the ideal $I$, and hence every entry of the $k^{th}$ power belongs to the ideal $I^{\lfloor k/m\rfloor}$.
\qed

\subsection{Formulation of the main result} We state a technical result 
needed to make the final arguments of the proof of Proposition \ref{Proposition:MainResultAlgebraicityProof}.
\subsubsection{Variables and rings}
Throughout the remainder of \S\ref{section:Gizmo} we fix a positive integer $n$
and we work with the family of independent (commuting) algebraic variables
$$\{u_i\}_{i=1}^n\cup\{v_i\}_{i=1}^{2n}\cup\{t,x,y\}.$$
Let 
$$u=(u_1,\dots,u_n)\;\;\mbox{and}\;\;v=(v_1,\dots,v_{2n}).$$
We write 
\begin{eqnarray*}
&&\CC[u]=\CC[u_1,\dots,u_n],\;
\CC[[u]]=\CC[[u_1,\dots,u_n]],\;\\
&&\CC[u,v]=\CC[u_1,\dots,u_n,v_1,\dots,v_{2n}],
\end{eqnarray*}
and so on. We use similar notation below for building up rings from the given variables without further comment.
Given, for example $P(u,v)\in \CC[u,v]^{2n}$, we denote by $\frac{\partial P}{\partial u}(u,v)$
the $2n$-by-$n$ matrix with entries $\frac{\partial P_i}{\partial u_j}(u,v)$. We use similar notation
for derivatives of vector functions below without further comment.

\subsubsection{Specialized matrix notation}
Given  $A\in \Mat_n(\CC((1/z)))$,
let 
\begin{eqnarray*}
A^\flat&=&\left[\begin{array}{ccccccc}
A(1,1)&\dots&A(n,1)&\dots &A(1,n)&\dots&A(n,n)
\end{array}\right]^\T\\
&\in&\Mat_{n^2\times 1}(\CC((1/z))).
\end{eqnarray*}
Note that
\begin{equation}\label{equation:flat}
(BA)^\flat=(I_n\otimes B)A^\flat
\end{equation}
for $B\in \Mat_n(\CC((1/z)))$.

 \subsubsection{Setup for the main result}
We are given an algebraic power series
\begin{equation}\label{equation:Typicalf}
f(t)=\sum_{i=0}^\infty c_it^i\in \CC[[t]]\;\;(c_i\in \CC)
\end{equation} 
and a matrix 
$$A\in \Mat_n(\CC((1/z)))$$ of negative spectral valuation, i.e., a matrix
satisfying the conditions
$$\lim_{k\rightarrow\infty}\val\, A^k=-\infty\;\;\mbox{and}\;\;e(A)\in (1/z)\CC[[1/z]]^{2n}$$ 
which by Proposition \ref{Proposition:NegativeSpectralValuation} are equivalent.
\begin{Proposition}\label{Proposition:Gizmo}
Notation and assumptions are as above.
For every $N\geq 0$ such that $\sum_{i=N}^\infty c_{i+N}t^i\in \CC[[t]]$ is nonsingular algebraic and divisible by $t^{2n}$,
there exist
$$\gamma\in (1/z)\CC[[1/z]]^{2n}\;\;\mbox{and}\;\;P(u,v)\in \CC[u,v]^{2n}$$
such that the following statements hold:
\begin{eqnarray}
\label{equation:Widget2}
&&P(e(A),\gamma)=0.\\
\label{equation:Widget3}
&&\det \frac{\partial P}{\partial v}(e(A),\gamma)\neq 0.\\
\label{equation:Widget4}
&&\left[\begin{array}{ccc}
(A^0)^\flat&\dots&(A^{2n-1})^\flat\end{array}\right]
\left(\frac{\partial P}{\partial v}(e(A),\gamma)\right)^{-1}
\frac{\partial P}{\partial u}(e(A),\gamma)=0.\\
\label{equation:Widget1}
&&\sum_{i=N}^\infty c_{i+N}\left[\begin{array}{cc}
A\otimes I_n&I_n\otimes I_n\\
0&I_n\otimes A
\end{array}\right]^i=
\sum_{i=1}^{2n}
\gamma_i\left[\begin{array}{cc}
A\otimes I_n&I_n\otimes I_n\\
0&I_n\otimes A
\end{array}\right]^{i-1}.
\end{eqnarray}
\end{Proposition}
\noindent The proof takes up the rest of \S\ref{section:Gizmo}
and is completed in \S\ref{subsection:GizmoFinish}. Note that every $N\geq 2n$
sufficiently large depending on $f(t)$ 
satisfies the hypotheses of 
Proposition \ref{Proposition:Gizmo}
by Proposition \ref{Proposition:Desingularization}.

\subsubsection{Remark}
In the application to the proof of Proposition \ref{Proposition:MainResultAlgebraicityProof} we will need to use several consequences
of the conclusions of Proposition \ref{Proposition:Gizmo} which are easy to check once written down
but might otherwise be obscure.
For the reader's convenience we write these consequences down. Firstly, we observe that the statement
\begin{eqnarray}
\label{equation:Widget4bis}
&&\left[\begin{array}{ccc}
(A^N)^\flat&\dots&(A^{N+2n-1})^\flat\end{array}\right]
\left(\frac{\partial P}{\partial v}(e(A),\gamma)\right)^{-1}
\frac{\partial P}{\partial u}(e(A),\gamma)=0
\end{eqnarray}
follows from statement \eqref{equation:Widget4} via statement  \eqref{equation:flat}. 
Secondly, we observe that \eqref{equation:Widget1}
implies
\begin{eqnarray}
\label{equation:Widget1.1}
&&\sum_{i=N}^\infty c_{i+N}\left[\begin{array}{cc}
A&B\\
0&A
\end{array}\right]^i=
\sum_{i=1}^{2n}
\gamma_i\left[\begin{array}{cc}
A&B\\
0&A
\end{array}\right]^{i-1}\;\;\mbox{for $B\in \Mat_n(\CC((1/z)))$}
\end{eqnarray}
by (so to speak) substituting $B$ for $\otimes$.
Finally, we observe that
\begin{eqnarray}
\label{equation:Widget1.3}
&&\sum_{i=2N}^\infty c_{i}A^i=
\sum_{i=1}^{2n}
\gamma_i A^{N+i-1}\;\;\mbox{and}\\\label{equation:Widget1.4}
&&\sum_{i=2N}^\infty \sum_{\nu=0}^{i-1}c_{i}A^\nu BA^{i-1-\nu}=
\sum_{i=1}^{2n}\sum_{\nu=0}^{N+i-2}\gamma_i
A^\nu BA^{i-2-\nu}
\;\;\mbox{for $B\in \Mat_n(\CC((1/z)))$.}
\end{eqnarray}
These last two statements are obtained by right-multiplying 
statement \eqref{equation:Widget1.1} on both sides by $\left[\begin{array}{cc}
A&B\\
0&A\end{array}\right]^N$ and expanding the matrix powers.

\subsubsection{Reduction of the proof of Proposition \ref{Proposition:Gizmo}}
After replacing $f(t)$ by \linebreak $\sum_{i=N}^\infty c_{i+N}t^i$,
we may assume without loss of generality that $N=0$. Thus we are making two
further special assumptions concerning $f(t)$ which for the sake of
clarity and convenient reference we write out explicitly.
Firstly, we are assuming that
\begin{equation}
\label{equation:Divisibility}
f(t)\in t^{2n}\CC[[t]].
\end{equation}
Secondly, we are assuming that there exists
$$F(x,y)\in \CC[x,y]$$
such that
\begin{equation}\label{equation:CleanAndSmooth}
F(t,f(t))=0\;\;\mbox{and}\;\;\frac{\partial F}{\partial y}(0,0)\neq 0.
\end{equation}
Note that the formula
\begin{equation}
\label{equation:CleanAndSmoothBis}
F(0,0)=\frac{\partial F}{\partial x}(0,0)=\cdots=\frac{\partial^{2n-1}F}{\partial x^{2n-1}}(0,0)=0
\end{equation}
follows straightforwardly from \eqref{equation:Divisibility} and
\eqref{equation:CleanAndSmooth}.

\subsection{A candidate for $\gamma$}

\subsubsection{A special polynomial}
Let
\begin{equation}\label{equation:DutDef}
D(u,t)=t^n+\sum_{i=1}^n (-1)^iu_it^{n-i}\in\CC[u,t].
\end{equation}
Note that the left side of \eqref{equation:CayleyHamilton}
equals  $D(e(A),A)$. This is the motivation for the definition of $D(u,t)$.

\subsubsection{Construction of $\varphi(u)$}
Perform Weierstrass division 
by $D(u,t)^2$  to obtain an identity
\begin{equation}\label{equation:varphiDef}
f(t)=\sum_{i=1}^{2n} \varphi_i(u)t^{i-1}+Q_1(u,t)D(u,t)^2,
\end{equation}
where 
$$\varphi(u)=(\varphi_1(u),\dots,\varphi_{2n}(u))\in u\CC[[u]]^{2n}\;\;\mbox{and}\;\;Q_1(u,t)\in \CC[[u,t]].$$
Note that we indeed have $\varphi(0)=0$ as one verifies by substituting $u=0$
on both sides of \eqref{equation:varphiDef} and using hypothesis
\eqref{equation:Divisibility}.

\begin{Lemma}\label{Lemma:Gizmo}
We have
\begin{eqnarray}\label{equation:PreWidget1}
&&\sum_{j=0}^\infty c_j
\left[\begin{array}{cc}
A\otimes I_n&I_n\otimes I_n\\
0&I_n\otimes A
\end{array}\right]^j=
\sum_{j=1}^{2n}
\varphi_j(e(A))\left[\begin{array}{cc}
A\otimes I_n&I_n\otimes I_n\\
0&I_n\otimes A
\end{array}\right]^{j-1}\;\;\mbox{and}\\
\label{equation:PreWidget5}
&&\left[\begin{array}{ccc}
(A^0)^\flat &\dots &(A^{2n-1})^\flat\end{array}\right]
\frac{\partial \varphi}{\partial u}(e(A))=0.
\end{eqnarray}
\end{Lemma}
\noindent  Thus the reasonable candidate 
for $\gamma$ is $\varphi(e(A))$.
 \proof
Note that \eqref{equation:PreWidget5} can be rewritten
 \begin{equation}
\label{equation:PreWidget2}
\sum_{i=1}^{2n}\frac{\partial \varphi_i}{\partial u_j}(e(A))A^{i-1}=0\;\;
\mbox{for $j=1,\dots,n$.}
\end{equation}
Note also that by differentiation we deduce from \eqref{equation:varphiDef} that
\begin{eqnarray}\label{equation:varphiDefbis}
&&\sum_{i=1}^{2n} \frac{\partial\varphi_i}{\partial u_j}(u)t^{i-1}\\
\nonumber&=&-\left(\frac{\partial Q_1}{\partial u_j}(u,t)D(u,t)
+2Q_1(u,t)\frac{\partial D}{\partial u_j}(u,t)\right)D(u,t)
\;\;\mbox{for $j=1,\dots,n$.}
\end{eqnarray}
Now suppose temporarily that $f(t)\in \CC[t]$.
Then we have $\varphi(u)\in u\CC[u]^{2n}$
and $Q_1(u,t)\in \CC[u,t]$
since high school division
in this case gives the same result as Weierstrass division.
Substituting 
$$(u,t)=\left(e(A),\left[\begin{array}{cc}
A\otimes I_n&I_n\otimes I_n\\
0&I_n\otimes A
\end{array}\right]\right)$$ into \eqref{equation:varphiDef} and using the Cayley-Hamilton Theorem \eqref{equation:CayleyHamilton},
 we obtain \eqref{equation:PreWidget1}.
Substituting 
$(u,t)=(e(A),A)$ into \eqref{equation:varphiDefbis} and using the Cayley-Hamilton Theorem \eqref{equation:CayleyHamilton}
  again,
we obtain \eqref{equation:PreWidget2}.
The general case follows by a routine approximation argument
based on Lemma \ref{Lemma:Approximation} and the remark immediately following.
\qed

\subsection{A candidate for $P(u,v)$}\subsubsection{Construction of the candidate}

Perform high school division of
$$F\left(t,\sum_{i=1}^{2n} v_it^{i-1}\right)\in \CC[u,v,t]$$
by $D(u,t)^2$ to obtain an identity
\begin{equation}\label{equation:PhiDef}
F\left(t,\sum_{i=1}^{2n} v_it^{i-1}\right)=\sum_{i=1}^{2n} P_i(u,v)t^{i-1}+Q_2(u,v,t)D(u,t)^2
\end{equation}
where
$$Q_2(u,v,t)\in \CC[u,v,t]$$
and
$$
P(u,v)=(P_1(u,v),\dots,P_{2n}(u,v))\in \CC[u,v]^{2n}.
$$
The latter is our candidate for $P(u,v)$.

\begin{Lemma}\label{Lemma:SuperGizmo}
Assumptions and notation are as above. Then the following statements hold:
\begin{eqnarray}
\label{equation:SuperGizmo3}
&&P(u,\varphi(u))=0.\\
\label{equation:SuperGizmo2}
&&\det \frac{\partial P}{\partial v}(0,0)\neq 0\;
\mbox{and hence}\;\frac{\partial P}{\partial v}(u,\varphi(u))\in \GL_{2n}(\CC[[u]]).\\
\label{equation:SuperGizmo4}
&&\frac{\partial \varphi}{\partial u}(u)=\left(
\frac{\partial P}{\partial v}(u,\varphi(u))
\right)^{-1}
\frac{\partial P}{\partial u}(u,\varphi(u)).
\end{eqnarray}
\end{Lemma}

\subsubsection{Proof of \eqref{equation:SuperGizmo3}}

Perform Weierstrass division of $F(x,y)$ by $y-f(x)$ in the
formal power series ring $\CC[[x,y]]$
to obtain the identity 
\begin{equation}
\label{equation:Weier1}
F(x,y)=(y-f(x))U(x,y)
\end{equation}
for some $U(x,y)\in \CC[[x,y]]$.
{\em A priori} one should add a remainder term $r(x)\in \CC[[x]]$ to the right side
but substitution of $y=f(x)$ on both sides
and the hypothesis \eqref{equation:CleanAndSmooth} shows that $r(x)=0$.
We then have the following chain of equalities:
\begin{eqnarray*}
&&\sum_{i=1}^{2n}P_i(u,\varphi(u))t^{i-1}+Q_2(u, \varphi(u),t)D(u,t)^2\;=\;F\left(t,\sum_{i=1}^{2n}\varphi(u)t^{i-1}\right)\\
&=&\left(\sum_{i=1}^{2n}\varphi(u)t^{i-1}-f(t)\right)U\left(t,\sum_{i=1}^{2n}\varphi_{i}(u)t^{i-1}\right)\\
&=&-Q_1(u,t)D(u,t)^2U\left(t,\sum_{i=1}^{n}\varphi(u)t^{i-1}\right).
\end{eqnarray*}
Justifications for the steps are as follows.
The first equality we obtain by substituting $v=\varphi(u)$ into \eqref{equation:PhiDef}.
The second equality we obtain by substituting into the factorization given in \eqref{equation:Weier1}.
The third equality we obtain by rearrangement of \eqref{equation:varphiDef}.
The equality between the extreme terms of the chain of equalities above forces \eqref{equation:SuperGizmo3} to hold
by the uniqueness of the remainder produced by Weierstrass division.
 \subsubsection{Proof of  \eqref{equation:SuperGizmo2}}
 Differentiation on both sides of \eqref{equation:PhiDef} with respect to $v_j$ followed by evaluation
at $u=v=0$ yields the relation
\begin{equation}\label{equation:WeierWildCard}
\frac{\partial F}{\partial y}\left(t,0\right)t^{j-1}= \sum_{i=1}^{2n}\frac{\partial P_i}{\partial v_j}(0,0)t^{i-1}+\frac{\partial Q_2}{\partial v_j}(0,0,t)t^{2n}.
\end{equation}
Now write
$$\frac{\partial F}{\partial y}(x,0)=\sum_{i=0}^\infty b_ix^i\;\;(b_i\in \CC).$$
By \eqref{equation:WeierWildCard} we have
$$\frac{\partial P_i}{\partial v_j}(0,0)=\left\{\begin{array}{rl}
b_{i-j}&\mbox{if $j\leq i$,}\\
0&\mbox{if $j>i$}
\end{array}\right.$$
for $i,j=1,\dots,2n$.
Thus we have
$$\det_{i,j=1}^{2n} \frac{\partial P_i}{\partial v_j}(0,0)=\left(\frac{\partial F}{\partial y}(0,0)\right)^{2n}.$$
The right side does not vanish by assumption \eqref{equation:CleanAndSmooth}.
Thus \eqref{equation:SuperGizmo2}  holds.

\subsubsection{Proof of \eqref{equation:SuperGizmo4}}
Formula \eqref{equation:SuperGizmo4} 
holds by implicit differentiation of formula \eqref{equation:SuperGizmo3}.
The proof of Lemma \ref{Lemma:SuperGizmo} is complete. 
\qed

\subsection{Completion of the proof of Proposition \ref{Proposition:Gizmo}}\label{subsection:GizmoFinish}
As noted above, we may assume that $N=0$
and hence that assumptions \eqref{equation:Divisibility},
\eqref{equation:CleanAndSmooth} and \eqref{equation:CleanAndSmoothBis} are in force.
Property \eqref{equation:Widget2} follows from 
formula \eqref{equation:SuperGizmo3}.
Property \eqref{equation:Widget3} follows
from formula \eqref{equation:SuperGizmo2}.
Property \eqref{equation:Widget4}
follows from formulas \eqref{equation:PreWidget5} and
\eqref{equation:SuperGizmo4}.
Property \eqref{equation:Widget1} follows 
from formula \eqref{equation:PreWidget1}.
The proof of Proposition \ref{Proposition:Gizmo} is complete.

 \qed

\section{Proof of the main result}\label{section:Endgame}
We finish the proof of Proposition \ref{Proposition:MainResultAlgebraicityProof} by checking hypotheses
in Proposition \ref{Proposition:AlgebraicityCriterion}, thereby completing the
proof of Theorem \ref{Theorem:MainResult}. 

\subsection{Review of the setup for Proposition \ref{Proposition:MainResultAlgebraicityProof}}
Let us start simply by repeating statements \eqref{equation:GSDE2} and \eqref{equation:GSDE3}
here for the reader's convenience:
\begin{equation}\label{equation:GSDE2bis}
I_{n}+a^{(0)} g+\sum_{\theta=1}^q \sum_{j=2}^\infty \kappa_j^{(\theta)}(a^{(\theta)} g)^j=0.
\end{equation}
\begin{eqnarray}
\label{equation:GSDE3bis}
&&\mbox{The linear map}\\
\nonumber&&\left(h\mapsto a^{(0)} h+\sum_{\theta=1}^q\sum_{j=2}^\infty
\sum_{\nu=0}^{j-1}
\kappa_j^{(\theta)}
(a^{(\theta)}g)^\nu (a^{(\theta)}h)(a^{(\theta)}g)^{j-1-\nu}
\right)\\
\nonumber&&:\Mat_n(\CC((1/z)))\rightarrow \Mat_n(\CC((1/z)))
\;\mbox{is invertible.}
\end{eqnarray}
Concerning the data appearing above, we have by \eqref{equation:Data} and \eqref{equation:GSDE5.5} 
that
\begin{equation}\label{equation:CoefficientCensus}
a^{(0)}\in \Mat_n(\CC(z)),\;\;a^{(1)},\dots,a^{(q)}\in \Mat_n(\CC)\;\;\mbox{and}\;\;
g\in \Mat_n(\CC((1/z))).
\end{equation}

\subsection{Application of Proposition \ref{Proposition:Gizmo}}

By \eqref{equation:Data}  and \eqref{equation:GSDE4}  we have that
$$
\mbox{$\sum_{j=2}^\infty \kappa_j^{(\theta)}t^j\in \CC[[t]]$ is algebraic
for $\theta=1,\dots,q$.}
$$
By assumption \eqref{equation:GSDE1} and Proposition \ref{Proposition:NegativeSpectralValuation} 
we have that
$$
\lim_{k\rightarrow\infty} \max_{\theta=1}^q \val\, (a^{(\theta)}g)^k=0\;\;\mbox{and}\;\;
e(a^{(1)}g),\dots,e(a^{(q)}g)\in (1/z)\CC[[1/z]]^n.
$$
Thus by Propositions \ref{Proposition:Desingularization}
and \ref{Proposition:Gizmo} along with the remarks immediately following the latter, there exist an integer $N\geq 2$
along with
$$\gamma^{(1)},\dots,\gamma^{(q)}\in (1/z)\CC[[1/z]]^{2n}\;\;\mbox{and}\;\;
P^{(1)}(u,v),\dots,P^{(q)}(u,v)\in \CC[u,v]^{2n}$$
such that 
\begin{eqnarray}\label{equation:Widget1Nought}
&&P^{(\theta)}(e(a^{(\theta)}g),\gamma^{(\theta)})=0,\\
\label{equation:Widget1NoughtBis}
&&\det\frac{\partial P^{(\theta)}}{\partial v}(e(a^{(\theta)}g),\gamma^{(\theta)})\neq 0,\\
\label{equation:Widget2Nought}
&&\sum_{j=2N}^\infty
\kappa_j^{(\theta)}(a^{(\theta)}g)^j=\sum_{j=1}^{2n}\gamma^{(\theta)}_j(a^{(\theta)}g)^{N+j-1},\\
\label{equation:Widget4Nought}
&&\left[\begin{array}{ccc}
((a^{(\theta)}g)^N)^\flat&\dots&((a^{(\theta)}g)^{N+2n-1})^\flat\end{array}\right]\\
\nonumber&&\times
\left(\frac{\partial P^{(\theta)}}{\partial v}(e(a^{(\theta)}g),\gamma^{(\theta)})\right)^{-1}
\frac{\partial P^{(\theta)}}{\partial u}(e(a^{(\theta)}g),\gamma^{(\theta)})=0,
\;\;\mbox{and}\\
\label{equation:Widget3Nought}
&&\sum_{j=2N}^\infty\sum_{\nu=0}^{j-1}
\kappa_j^{(\theta)}(a^{(\theta)}g)^\nu(a^{(\theta)}h)(a^{(\theta)}g)^{j-1-\nu}\\
\nonumber&&=\sum_{j=1}^{2n}\sum_{\nu=0}^{N+j-2}
\gamma_j^{(\theta)}(a^{(\theta)}g)^\nu(a^{(\theta)}h)(a^{(\theta)}g)^{N+j-2-\nu}
\end{eqnarray}
for $\theta=1,\dots,q$ and any $h\in \Mat_n(\CC((1/z)))$.

\subsection{Polynomial version of \eqref{equation:GSDE2bis}}
We embed  \eqref {equation:GSDE2bis} into a system of $3nq+n^2$ polynomial equations in $3nq+n^2$ variables with coefficients in $\CC(z)$.
 \subsubsection{Variables}
 We employ the family of variables
$$\{U_i\}_{i=1}^{qn}\cup\{V_i\}_{i=1}^{2qn}\cup\left\{\xi_i\right\}_{i=1}^{n^2}.$$
Let
$$U=(U_1,\dots,U_{qn}),
\;\;V=(V_1,\dots,V_{2qn}),\;\;\mbox{and}\;\;
\xi=(\xi_1,\dots,\xi_{n^2}).
$$
Let
$$\Xi=\left[\begin{array}{ccccc}
\xi_1&\dots&\xi_{n^2-n+1}\\
\vdots&&\vdots\\
\xi_n&\dots&\xi_{n^2}
\end{array}\right]\in \Mat_n(\CC[\xi]).
$$
We break the $U$'s and $V$'s down into groups by introducing the following notation:
\begin{eqnarray*}
U^{(\theta)}_i&=&U_{i+(\theta-1)q}\;\;\mbox{for $i=1,\dots,n$ and $\theta=1,\dots,q$.}\\
V^{(\theta)}_i&=&V_{i+(\theta-1)q}\;\;\mbox{for $i=1,\dots,2n$ and $\theta=1,\dots,q$.}\\
U^{(\theta)}&=&(U_1^{(\theta)},\dots,U_n^{(\theta)})\;\;\mbox{and}\;\;V^{(\theta)}\;=\;
(V_1^{(\theta)},\dots,V_{2n}^{(\theta)})\;\;\mbox{for $\theta=1,\dots,q$.}
\end{eqnarray*}
\subsubsection{A special matrix}
We define
\begin{eqnarray*}
&&\HHH(V,\xi)\\
&=&I_n+a^{(0)}\Xi+\sum_{\theta=1}^q\left(\sum_{k=2}^{2N-1}
\kappa_k^{(\theta)}(a^{(\theta)}\Xi)^k+
\sum_{k=1}^{2n}V_k^{(\theta)}(a^{(\theta)}\Xi)^{N+k-1}\right)\\
&\in& \Mat_n(\CC(z)[V,\xi]).
\end{eqnarray*}
\subsubsection{Polynomials}
We define $3qn+n^2$ polynomials belonging to $\CC(z)[U,V,\xi]$ as follows.
\begin{eqnarray*}
&&F_{i+n(\theta-1)}(U,\xi)\\
&=&F^{(\theta)}_i(U^{(\theta)},\xi)=U_i^{(\theta)}-e_i(a^{(\theta)}\Xi)\;\;\mbox{for $i=1,\dots,n$ and $\theta=1,\dots,q$,}\\\\
&&G_{i+2n(\theta-1)}(U,V)\\
&=&P^{(\theta)}(U^{(\theta)},V^{(\theta)})\;\;\mbox{for $i=1,\dots,2n$ and $\theta=1,\dots,q$, and}\\\\
&&H_{i+n(j-1)}(V,\Xi)\\
&=&\HHH(V,\Xi)(i,j)\;\;\mbox{for $i,j=1,\dots,n$.}
\end{eqnarray*}
\subsubsection{Presentation of the system of equations}
Let
\begin{eqnarray*}
F(U,\xi)&=&\left[\begin{array}{ccc}F_1(U,\xi)&\dots&F_{qn}(U,\xi)\end{array}\right]^\T,\\
G(U,V)&=&\left[\begin{array}{ccc}G_1(U,V)&\dots&G_{2qn}(U,V)\end{array}\right]^\T,\\
H(V,\xi)&=&\left[\begin{array}{ccc}H_1(V,\xi)&\dots&H_{n^2}(V,\xi)\end{array}\right]^\T.
\end{eqnarray*}
Then our system of polynomial equations takes the form
\begin{equation}\label{equation:TheSystem}
F(U,\xi)=0,\;\;G(U,V)=0,\;\;H(V,\xi)=0.
\end{equation}
Note that this system has all coefficients in $\CC(z)$
by \eqref{equation:CoefficientCensus} and the definitions.

\subsubsection{The solution $\Upsilon_0$}
We claim that the following formulas specify a solution over $\CC((1/z))$ 
of the system of equations \eqref{equation:TheSystem}:
\begin{eqnarray}\label{equation:Group1Nought}
&&U_i^{(\theta)}=e_i(a^{(\theta)}g)\;\;
\mbox{for $i=1,\dots,n$ and $\theta=1,\dots,q$.}\\
\label{equation:Group2Nought}
&&V_i^{(\theta)}=\gamma_i^{(\theta)}\;\;\mbox{for $i=1,\dots,2n$ and $\theta=1,\dots,q$.}\\
\label{equation:Group3Nought}
&&\xi_{i+n(j-1)}=\Xi(i,j)=g(i,j)\;\;\mbox{for $i,j=1,\dots,n$.}
\end{eqnarray}
The equation $F(U,\xi)=0$ is obviously satisfied.
The equation $G(U,V)=0$ is satisfied because it merely restates
the system of equations \eqref{equation:Widget1Nought}. Finally, one verifies
that $H(U,\xi)=0$ is satisfied by using \eqref{equation:GSDE2bis},
 \eqref{equation:Widget2Nought},  \eqref{equation:Group1Nought}, and  \eqref{equation:Group2Nought}.
 The claim is proved.
 The solution of  \eqref{equation:TheSystem} specified
 by \eqref{equation:Group1Nought}, \eqref{equation:Group2Nought}, and \eqref{equation:Group3Nought}
   will be denoted by $\Upsilon_0$.
\subsection{Analysis of the Jacobian determinant }
Now we study the Jacobian matrix
\begin{equation}\label{equation:BigMatrix}
\left[\begin{array}{ccc}
\displaystyle\frac{\partial F}{\partial U}(U,\xi)&0&\displaystyle\frac{\partial F}{\partial \xi}(U,\xi)\\\\
\displaystyle\frac{\partial G}{\partial U}(U,V)&\displaystyle\frac{\partial G}{\partial V}(U,V)&0\\\\
0&\displaystyle\frac{\partial H}{\partial V}(V,\xi)&\displaystyle\frac{\partial H}{\partial \xi}(V,\xi)
\end{array}\right]\in \Mat_{n^2+3qn}(\CC(z)[U,V,\xi])
\end{equation}
for the system of equations \eqref{equation:TheSystem}.
Let
\begin{equation}\label{equation:BigMatrixb}
\left[\begin{array}{ccc}
I_n&0&b_{13}\\
b_{21}&b_{22}&0\\
0&b_{32}&b_{33}
\end{array}\right]\in \Mat_{n^2+3qn}(\CC((1/z)))
\end{equation}
be the result of evaluating \eqref{equation:BigMatrix} at the point $\Upsilon_0$.
To prove Proposition \ref{Proposition:MainResultAlgebraicityProof}
and thereby to complete the proof of Theorem \ref{Theorem:MainResult},
we have by Proposition \ref{Proposition:AlgebraicityCriterion}
only to prove that the determinant of the matrix \eqref{equation:BigMatrixb} does not vanish.
Now provided that $\det b_{22}\neq 0$, we have a matrix identity
\begin{eqnarray*}
&&\left[\begin{array}{ccc}
I_n&0&b_{13}\\
b_{21}&b_{22}&0\\
0&b_{32}&b_{33}
\end{array}\right]\left[\begin{array}{ccc}
I_n&0&0\\
-b_{22}^{-1}b_{21}&I_{2n}&0\\
0&0&I_{n^2}
\end{array}\right]\left[\begin{array}{ccc}
I_n&0&-b_{13}\\
0&I_{2n}&0\\
0&0&I_{n^2}\end{array}\right]\\
\nonumber&=&\left[\begin{array}{ccc}
I_n&0&0\\
0&b_{22}&0\\
-b_{32}b_{22}^{-1}b_{21}&b_{32}&b_{32}b_{22}^{-1}b_{21}b_{13}+b_{33}\end{array}\right].
\end{eqnarray*}
Thus it will be enough to prove that
\begin{eqnarray}\label{equation:Clincher1}
&&\det b_{22}\neq 0,\\
\label{equation:Clincher2}
&&b_{32}b_{22}^{-1}b_{21}=0,\;\;\mbox{and}\\
\label{equation:Clincher3}
&&\det b_{33}\neq 0.
\end{eqnarray}

\subsection{Proof of \eqref{equation:Clincher1}}
We have by the definitions
\begin{equation}\label{equation:b22}
b_{22}=\sum_{\theta=1}^q
\ebold_{\theta\theta}\otimes\frac{\partial P^{(\theta)}}{\partial v}(e(a^{(\theta)}g),\gamma^{(\theta)}).\end{equation}
Thus \eqref{equation:Clincher1} holds by \eqref{equation:Widget1NoughtBis}.

\subsection{Proof of \eqref{equation:Clincher2}}
For $b_{21}$ and $b_{32}$ we have formulas similar to \eqref{equation:b22}, namely
\begin{eqnarray*}
b_{21}&=&\sum_{\theta=1}^q
\ebold_{\theta\theta}\otimes
\frac{\partial P^{(\theta)}}{\partial u}(e(a^{(\theta)}g),\gamma^{(\theta)})\;\;\mbox{and}\\
b_{32}&=&\sum_{\theta=1}^q
\ebold_{\theta\theta}\otimes\left[\begin{array}{ccc}
((a^{(\theta)}g)^N)^\flat & \dots & ((a^{(\theta)}g)^{N+2n-1})^\flat\end{array}\right].
\end{eqnarray*}
Thus \eqref{equation:Clincher2} holds by
\eqref{equation:Widget4Nought}.

\subsection{Proof of \eqref{equation:Clincher3}}
We have for $i,j=1,\dots,n$ that 
\begin{eqnarray*}
&&\frac{\partial \HHH(V,\xi)}{\partial \xi_{i+(j-1)n}}\\
&=&a^{(0)}\ebold_{ij}+\sum_{\theta=1}^q\sum_{k=2}^{2N-1}
\sum_{\nu=0}^{k-1}\kappa_k^{(\theta)}(a^{(\theta)}\Xi)^\nu(a^{(\theta)}\ebold_{ij})(a^{(\theta)}\Xi)^{k-1-\nu}\\
\nonumber&&+\sum_{\theta=1}^q
\sum_{k=1}^{2n}\sum_{\nu=0}^{N+k-2}V_k^{(\theta)}(a^{(\theta)}\Xi)^\nu
(a^{(\theta)}\ebold_{ij})(a^{(\theta)}\Xi)^{N+k-2-\nu}
\end{eqnarray*}
and hence after evaluating both sides at $\Upsilon_0$ and using \eqref{equation:Widget3Nought},
we find that $b_{33}$ is a matrix describing with respect to the basis
$$\ebold_{11},\dots,\ebold_{n1},\dots,\ebold_{1n},\dots,\ebold_{nn}\in\Mat_n(\CC((1/z)))$$ the invertible linear map considered in \eqref{equation:GSDE3bis}.
Thus \eqref{equation:Clincher3} holds. Thus the proof of Proposition \ref{Proposition:MainResultAlgebraicityProof} is complete and with it the proof of Theorem \ref{Theorem:MainResult}.
\\

\textbf{Acknowledgements:}  I  thank Serban Belinschi,
 J. William Helton, Tobias Mai,
and Roland Speicher for communications concerning the self-adjoint linearization trick
and symmetric realization.
I  thank Steven Lalley for communications concerning algebraicity in relation to random walk. 
 I thank Christine Berkesch Zamaere for communications concerning
algebraicity criteria and commutative algebra.

\end{document}